\documentclass[11pt]{amsart}
\usepackage{amsmath}
\usepackage{amssymb}
\usepackage{amsfonts}
\usepackage{amsthm}
\usepackage{enumerate}
\usepackage[mathscr]{eucal}
\usepackage{verbatim}
\usepackage{amsthm}
\usepackage{amscd}
\usepackage[mathscr]{eucal}
\usepackage{appendix}
\usepackage{tikz}
\usepackage{hyperref}
\usepackage{cleveref}

\numberwithin{equation}{section}

\newtheorem{innercustomgeneric}{\customgenericname}
\providecommand{\customgenericname}{}

\newcommand{\newcustomtheorem}[2]{\newenvironment{#1}[1]
  {\renewcommand\customgenericname{#2}
   \renewcommand\theinnercustomgeneric{##1}\innercustomgeneric}{\endinnercustomgeneric}}

\newcustomtheorem{customthm}{Theorem}

\newcommand{\newcustomlemma}[2]{\newenvironment{#1}[1]
  {\renewcommand\customgenericname{#2}
   \renewcommand\theinnercustomgeneric{##1} \innercustomgeneric}{\endinnercustomgeneric}}

\newcustomlemma{customlemma}{Lemma}

\newcustomlemma{customproposition}{Proposition}

\newcustomlemma{customclaim}{Claim}

\theoremstyle{plain}
\newtheorem{theorem}{Theorem}

\newtheorem{lemma}[theorem]{Lemma}
\newtheorem{proposition}[theorem]{Proposition}
\newtheorem{remark}{Remark}

\newtheorem*{theorem*}{Theorem}
\newtheorem*{lemma*}{Lemma}
\newtheorem*{proposition*}{Proposition}
\newtheorem*{corollary*}{Corollary}
\newtheorem*{remark*}{Remark} 
\newtheorem*{remarks*}{Remarks}
\newtheorem*{conj*}{Conjecture}

\def\S{\mathbb{S}}

\def\R{{\mathbb R}}
\def\Z{{\mathbb Z}}

\def\B{{\mathbb B}}

\newcommand{\LL}{\mathcal{L}}

\newcommand{\WW}{\mathcal{W}}
\newcommand{\TT}{\mathcal{T}}

\newcommand{\bbz}{\mathbb{Z}}

\newcommand{\bbs}{\mathbb S}

\newcommand{\bbr}{\mathbb{R}}

\newcommand{\bbn}{\mathbb{N}}

\newcommand{\kkk}{\vec{\boldsymbol{k}}}
\newcommand{\xxxi}{\vec{\boldsymbol{\xi}\;}}
\newcommand{\xxx}{\vec{\boldsymbol{x}}}
\newcommand{\yyy}{\vec{\boldsymbol{y}}}
\newcommand{\zzz}{\vec{\boldsymbol{z}}}

\newcommand{\GGG}{\vec{\boldsymbol{G}}}
\def\000{\vec{\boldsymbol{0}}}

\newcommand{\II}{\mathcal{I}}
\newcommand{\UU}{\mathcal{U}}

\def\ga{\gamma}

\def\s{\sigma}

\def\ep{\epsilon}
\def\Om{\Omega}
\newcommand{\la}{\lambda}

\newcommand{\q}{\quad}
\newcommand{\qq}{\qquad}

\DeclareFontFamily{U}{mathx}{\hyphenchar\font45}
\DeclareFontShape{U}{mathx}{m}{n}{
	<5> <6> <7> <8> <9> <10>
	<10.95> <12> <14.4> <17.28> <20.74> <24.88>
	mathx10
}{}

\def\wh{\widehat}

\newcommand{\wt}{\widetilde}

\newcommand{\supp}{\mathrm{supp}}

\makeatletter
\@namedef{subjclassname@2020}{\textup{2020} Mathematics Subject Classification}
\makeatother

\allowdisplaybreaks

\makeindex         

\begin{document}

\author{Stefanos Lappas}
\address{S. Lappas, Department of Mathematical Analysis, Faculty of Mathematics and
Physics, Charles University, Sokolovsk\'a 83, 186 75 Praha 8, Czech Republic}
\email{stefanos.lappas@matfyz.cuni.cz}

\author{Bae Jun Park}
\address{B. Park, Department of Mathematics, Sungkyunkwan University, Suwon 16419, Republic of Korea}
\email{bpark43@skku.edu}

\thanks{S. Lappas was supported by the Primus research programme PRIMUS/21/SCI/002 of Charles University and by the Foundation for Education and European Culture, founded by Nicos and Lydia Tricha. B. Park was supported by the National Research Foundation (Republic of
Korea) grant RS-2025-20512969 and by the Open KIAS Center at Korea Institute for Advanced Study. }

\title[Sharp bilinear estimates for maximal singular integrals]{Sharp bilinear estimates for maximal singular integrals with kernels in weighted $L^q$ spaces} 


\begin{abstract}
In this paper, we study the boundedness properties of the (dyadic) maximal bilinear operator associated with rough homogeneous kernels on $\mathbb{R}$.
We establish sharp $L^{p_1}(\mathbb{R}) \times L^{p_2}(\mathbb{R}) \to L^{p}(\mathbb{R})$ estimates in the full quasi-Banach range of exponents $1 < p_1, p_2 < \infty$ and $1/2 < p < \infty$.
Our approach extends and unifies several recent contributions, including those of Honz\'ik, the first author, and Slav\'ikova \cite{Ho_La_Sl_submitted}, as well as the second author \cite{Park2025} in the bilinear and in the one-dimensional settings, by allowing the angular component $\Omega$ of the kernel to belong to weighted $L^q$-spaces on $\mathbb{S}^1$.
\end{abstract}

\maketitle

\section{Introduction}

The work of Calder\'on and Zygmund \cite{CZ52} introduced singular integral operators of the  form
\begin{equation}\label{eq:linear}
  \LL f(x):=\text{p.v.}\int_{\R^n}K(y)f(x-y)~dy.
\end{equation}
Since then, numerous remarkable boundedness results of the operator $\LL$ appeared in the literature. In particular, if the kernel of \eqref{eq:linear} takes the form $K(y)=\widetilde\Omega(\frac{y}{|y|})|y|^{-n}$, where the function $\widetilde\Omega$ has vanishing integral over the unit sphere $\S^{n-1}$ and, in addition, belongs to $L\log L(\S^{n-1})$, then the same authors of \cite{CZ52} obtained in \cite{CZ56} the boundedness of $\LL$ on $L^p(\R^n)$ for $1<p<\infty$. 
Subsequent developments refined this result through the use of Hardy space theory, notably in \cite{CW77,C79}, and further extensions, including endpoint estimates for $\LL$, were later established in \cite{C88, CRdF88, Ho1988, Se1996}.

As explained below, we focus on establishing estimates for maximal bi-sublinear operators associated with the bilinear version of operator \eqref{eq:linear} in the one-dimensional setting $n=1$. To this end, we first recall the notion of the bilinear rough singular operators. Suppose that $\Omega$ is an integrable function on the unit circle $\mathbb{S}^{1}$ with the mean value zero property 
\begin{equation}\label{vanishingmtcondition}
  \int_{\mathbb{S}^{1}}\Omega(\yyy')~ d\sigma(\yyy')=0
\end{equation} 
where $d\sigma$ stands for the circular measure on $\mathbb{S}^{1}$, $\yyy:=(y_1,y_2)\in \bbr^2$ and $\yyy':=\frac{\yyy}{|\yyy|}\in \mathbb{S}^{1}$.
We define the corresponding truncated bilinear operator $\LL_{\Om}^{(\epsilon)}$  by
\begin{equation}\label{eq:truncated}
  \mathcal L^{(\ep)}_{\Om}\big(f_1,f_2\big)(x):=\int_{|\yyy|>\epsilon}{K(\yyy)f_1(x-y_1)f_2(x-y_2)}~d\yyy, \qquad x\in\bbr
\end{equation}
where 
\begin{equation*}
  K(\yyy):=\frac{\Omega(\yyy')}{|\yyy|^{2}}, \qquad \yyy \neq \000.     
\end{equation*}
By taking the limit $\epsilon\searrow0$ in \eqref{eq:truncated}, we can define the bilinear rough singular integral operator $\LL_{\Om}$ as follows:
\begin{equation}\label{eq:rough}
  \LL_{\Om}(f_1,f_2)(x):=\lim_{\epsilon\searrow0}\mathcal L^{(\ep)}_{\Om}\big(f_1,f_2\big)(x)=\text{p.v.}\int_{\R^2}{K(\yyy)f_1(x-y_1)f_2(x-y_2)}~d\yyy, \qquad x\in\bbr.    
\end{equation}
Notice that \eqref{eq:truncated} and \eqref{eq:rough} are well-defined for any Schwartz functions $f_1,f_2$ on $\bbr$.

When $\Omega$ is a function of bounded variation on $\S^1$, then the first boundedness property for the bilinear singular integral operator $\LL_{\Om}$ was established by Coifman and Meyer  \cite{CoMe1975}. Extensions of this result were subsequently obtained by Kenig and Stein  \cite{KS99} and by Grafakos and Torres \cite{GT02}. 

According to \cite[pages 162 -- 163]{GT02}, if $\Omega$ is merely integrable on the unit circle $\S^1$, but is odd, the bilinear rough singular integral operator $\LL_{\Om}$ is intimately connected with the celebrated (directional) bilinear Hilbert transform whose  $L^{p_1}(\R) \times L^{p_2}(\R) \rightarrow L^{p}(\R)$ boundedness and its uniform estimates were established in the influential works of Lacey--Thiele \cite{LT97,LT99} and \cite{GL, Li, Thi02, UW}, respectively.

Grafakos, He and Honz{\'i}k \cite{Gr_He_Ho2018} proved the $L^{p_1}(\R) \times L^{p_2}(\R) \rightarrow L^{p}(\R)$ boundedness of $\LL_{\Om}$ when $1<p_1, p_2<\infty$, $1/2<p<\infty$, $1/p=1/p_1+1/p_2$ and $\Omega\in L^\infty(\S^1)$. In the same work a proof for the ``local $L^2$'' boundedness of $\LL_{\Om}$, namely $2\leq p_1, p_2\leq\infty$, $1\leq p\leq 2$, was  provided for $\Omega\in L^2(\S^{1})$. A further improvement of this result to all $\Omega\in L^q(\S^1)$ for $q>4/3$ was achieved in \cite{Gr_He_Sl2020} and this result was further extended to exponents $1<p_1,p_2<\infty$ and $1/2<p<\infty$, by He and the second author \cite{He_Park2023} when $\Omega\in L^{q}(\S^1)$ for $q>\max\big(\frac{4}{3},\frac{p}{2p-1}\big)$. Multilinear extensions of the latter result appeared in the series of works \cite{Gr_He_Ho_Park2023, Gr_He_Ho_Park_JLMS, Do_Sl2024}. We point out that the work of Dosidis and Slav\'ikov\'a \cite{Do_Sl2024} obtained that  
$\LL_{\Om}$ is $L^{p_1}(\R) \times L^{p_2}(\R) \rightarrow L^{p}(\R)$ bounded for all $\Omega\in L^q(\S^1)$ if and only if $1/p+1/q<2$, where $1<p_1, p_2<\infty$, $1/2<p<\infty$, $1/p=1/p_1+1/p_2$, and $q>1$. 
An alternative proof of this result, based on a key trilinear estimate with decay in the frequency parameter for a bilinear local operator acting on functions that are Fourier localized to an annulus, was also given in \cite[Theorem 1.1]{BS2025}. It is also worth mentioning that further refinements of the boundedness results for the operator $\LL_{\Om}$, allowing $\Omega$ to belong to a larger class of Orlicz spaces, were established in \cite[Theorem 3]{Do_Park_Sl_submitted} and \cite[Theorem 1.3]{BS2025}.

Very recently, Honz\'ik, the first author, and Slav\'ikov\'a \cite[Theorem 1]{Ho_La_Sl_submitted} showed that $\LL_{\Om}$ is $L^{p_1}(\R) \times L^{p_2}(\R) \rightarrow L^{p}(\R)$ bounded in the optimal quasi-Banach range of exponents $1<p_1,p_2<\infty$ and $1/2<p<\infty$, provided that $\Omega$ is supported away from the diagonal $\{(\theta_1,\theta_2)\in\S^1:\theta_1=\theta_2\}$. Furthermore, \cite[Theorem 1.4]{BS2025} offered an alternative proof of \cite[Theorem 1]{Ho_La_Sl_submitted} with minor  modifications to the argument used in the aforementioned result \cite[Theorem 1.1]{BS2025}. In fact, a more general result was obtained in
\cite[Theorem 2]{Ho_La_Sl_submitted} where the assumption $\Omega\in L^q(\mathbb{S}^1,u_A^q)$ for some $q>1$ was imposed. 
In order to describe this setting precisely, we introduce some notation. For $A>0$, let
\begin{equation*}
u_A(\theta_1,\theta_2):=\frac{1}{|\theta_1-\theta_2|^{A}} \quad \text{for} \quad (\theta_1,\theta_2)\in \S^1,
\end{equation*}
and then we define
\begin{equation*}
\Vert \Omega\Vert_{L^q(\mathbb{S}^1,u_A^q)}:=\bigg(\int_{\mathbb{S}^1}\big| \Omega(\theta_1,\theta_2)\big|^q u_A(\theta_1,\theta_2)^q\; d\sigma(\theta_1,\theta_2) \bigg)^{\frac{1}{q}}.
\end{equation*}

\begin{customthm}{A}\cite[Theorem 2]{Ho_La_Sl_submitted}\label{thm:submitted2024}
Let $1<p_1,p_2<\infty$, $\frac{1}{2}<p<\infty$ with $\frac{1}{p}=\frac{1}{p_1}+\frac{1}{p_2}$, and let $q>1$ be such that $\frac{1}{p}+\frac{1}{q}\geq 2$ and $\Omega\in L^{q}(\S^{1},u_A^q)$ with \eqref{vanishingmtcondition}. 
Then
\begin{equation*}
  A>\frac{1}{p}+\frac{1}{q}-2
\end{equation*}
if and only if there exists a constant $C=C(p_1,p_2,q,A)$ such that
\begin{equation*}
  \big\Vert \LL_{\Omega}(f_1,f_2)\big\Vert_{L^{p}(\R)}\leq C\Vert\Omega\Vert_{L^q(\S^{1},u_A^q)}\Vert f_1\Vert_{L^{p_1}(\R)}\Vert f_2\Vert_{L^{p_2}(\R)}.
\end{equation*}
for Schwartz functions $f_1,$ $f_2$ on $\bbr$.
\end{customthm}

We proceed by recalling the definition of maximal bi-sublinear operator $\LL_{\Om}^{*}$, given by
\begin{equation*}
\LL_{\Om}^{*}\big(f_1,f_2\big)(x) := 
\sup_{\ep>0 } \big| \LL_{\Om}^{(\ep)}\big(f_1, f_2\big)(x)\big|, \qquad x\in\bbr
\end{equation*}
for Schwartz functions $f_1,f_2$ on $\bbr$. The initial study of bilinear estimates for the maximal operator $\LL_{\Om}^{*}$ was carried out by Buri\'ankov\'a and Honz\'ik \cite{BuHo2019}. For multilinear extensions of this result, we refer the reader to \cite{Gr_He_Ho_Park2024, Park2025}. Now, for our purposes we recall from \cite[Theorem 1]{Park2025} the following bilinear estimate for $\LL_{\Om}^{*}$ in the one dimensional setting.   
\begin{customthm}{B}\cite[Theorem 1]{Park2025}\label{thm:Park2025}
Let $1<p_1,p_2<\infty$ and $\frac{1}{2}<p<\infty$ with $\frac{1}{p}=\frac{1}{p_1}+\frac{1}{p_2}$. Suppose that $\Omega\in L^{q}(\mathbb{S}^{1})$, $q>1$, with \eqref{vanishingmtcondition}. If $\frac{1}{p}+\frac{1}{q}<2$, then there exists a constant $C=C(p_1,p_2,q)>0$ such that
\begin{equation*}
 \big\|\LL_{\Om}^{*}(f_1, f_2)\big\|_{L^{p}(\mathbb{R}) } \le C\Vert \Omega\Vert_{L^{q}(\mathbb{S}^{1})} \Vert f_1\Vert_{L^{p_1}(\bbr)}\Vert f_2\Vert_{L^{p_2}(\bbr)}
 \end{equation*}
for Schwartz functions $f_1,f_2$ on $\bbr$.
\end{customthm}

In this paper, we consider a dyadic version of the maximal bi-sublinear operators.
We define
\begin{equation*}
  \LL_{\Omega}^{\mathrm{dyad},*}\big(f_1,f_2 \big)(x):=\sup_{\rho\in\mathbb{Z}}\big| \LL_{\Om}^{(2^{\rho})}\big(f_1, f_2\big)(x)\big|, \qquad x\in\bbr
\end{equation*}
for Schwartz functions $f_1,f_2$ on $\bbr$. Then our main result reads as follows.

\begin{theorem}\label{thm:main result}
Let $1<p_1,p_2<\infty$, $\frac{1}{2}<p<\infty$ with $\frac{1}{p}=\frac{1}{p_1}+\frac{1}{p_2}$, and let $q>1$ be such that $\frac{1}{p}+\frac{1}{q}\geq 2$ and $\Omega\in L^{q}(\S^{1},u_A^q)$ with \eqref{vanishingmtcondition}. 
If
\begin{equation}\label{sharpcond}
  A>\frac{1}{p}+\frac{1}{q}-2,
\end{equation}
then  there exists a constant $C=C(p_1,p_2,q,A)$ such that
\begin{equation*}
  \big\Vert \LL_{\Omega}^{\mathrm{dyad},*}(f_1,f_2)\big\Vert_{L^{p}(\R)}\leq C\Vert\Omega\Vert_{L^q(\S^{1},u_A^q)}\Vert f_1\Vert_{L^{p_1}(\R)}\Vert f_2\Vert_{L^{p_2}(\R)}
\end{equation*}
for Schwartz functions $f_1, f_2$ on $\bbr$.
\end{theorem}

\begin{remark}
We note that condition \eqref{sharpcond} is sharp, as the same counter-example used in the proof of Theorem \ref{thm:submitted2024} can be applied. Clearly, $\LL_{\Omega}^{\mathrm{dyad},*}(f_1,f_2)$  pointwise dominates $\LL_{\Omega}(f_1,f_2)$ almost everywhere. This observation, combined with the embedding $L^q(\S^{1},u_A^q)\hookrightarrow L^{q}(\S^{1})$, implies that our Theorem \ref{thm:main result} extends the results of \cite[Theorems 1 and 2]{Ho_La_Sl_submitted} and \cite[Theorems 1 and 2]{Park2025}. 
\end{remark}

The proof of Theorem \ref{thm:main result} relies on a pointwise decomposition of the maximal operator  $\LL_{\Omega}^{\mathrm{dyad},*}$ into two auxiliary maximal ones $\mathcal{E}^*$ and $\mathcal{L}_{\Omega}^{\sharp}$ as in \eqref{eq:maindecomp.}.
The latter, $\mathcal{L}_{\Omega}^{\sharp}$, corresponds to a smoothly truncated version of the kernel, while  $\mathcal{E}^*$ captures the additional contribution that arises when the rough truncation $\chi_{|\yyy|>2^{\rho}}$ is replaced by its smooth one $1-\widehat{\Phi}(2^{-\rho}\yyy)$.
This decomposition enables us to analyze the boundedness properties of $\mathcal{E}^*$ and $\mathcal{L}_{\Omega}^{\sharp}$ separately, and it also provides a natural framework in which the dyadic kernel decomposition of Duoandikoetxea and Rubio de Francia \cite{Du_Ru1986} can be effectively employed.

In contrast to the earlier work \cite{Gr_He_Ho_Park2024}, which treated similar multi-sublinear maximal operators associated with rough singular integrals,
we encounter substantial new difficulties in the present setting. In the previous case, the term $\mathcal{E}^*$ could be handled rather directly under the assumption $\Omega\in L^q$ for sufficiently large $q>1$, by invoking \cite[Lemma 3.1]{Gr_He_Ho_Park2024}.
Here, however, we must deal with the more general situation $\Omega\in L^q(\mathbb{S}^1,u_A^q)$ for arbitrary $q>1$, where the previous argument is no longer applicable.

Developing a new approach to control $\mathcal{E}^*$ under this general setting constitutes one of the main contributions of the paper. This part of the analysis is organized into six key propositions (Propositions \ref{muge0tcpg1}-\ref{prop:intmaxE1}), whose proofs are postponed to later sections to maintain the continuity of exposition.

For the term $\mathcal{L}_{\Omega}^{\sharp}$, the dyadic kernel decomposition of \cite{Du_Ru1986} fits naturally into the framework, allowing us to follow the general strategy developed in \cite{Gr_He_Ho_Park2024}. To eliminate the dependence of the parameter $q$ on the exponent $p$, we establish an additional endpoint estimate (Proposition \ref{lommushweakre}) and then apply interpolation to cover the full range of boundedness.

Combining the estimates for both $\mathcal{E}^*$ and $\mathcal{L}_{\Omega}^{\sharp}$, we obtain the desired $L^{p_1}(\R)\times L^{p_2}(\R)\to L^p(\R)$ boundedness under the condition $A>\frac{1}{q}$. Finally, to further refine this, by regarding the bi-sublinear operator involving $\Omega$ as a tri-linearizable operator with $\Omega$ treated as a variable function, we apply the weighted interpolation theorem of Cao, Olivo, and Yabuta \cite{COY22} to refine the condition $A>\frac{1}{q}$ to the sharp range $A>\frac{1}{p}+\frac{1}{q}-2$, as stated in \eqref{sharpcond}. This final interpolation step is carried out in Subsection \ref{Sec: Section 3.3}.

\hfill

Our paper is organized as follows. Section \ref{Sec: Section 2} provides essential preliminaries, including key concepts from wavelet theory and interpolation that will be used later. In Section \ref{Sec: Section 3}, we present the proof of Theorem \ref{thm:main result}. After establishing a suitable pointwise decomposition of the maximal operator, we derive two main estimates for the auxiliary maximal functions introduced earlier and prove them through a series of propositions. The subsequent sections (Sections \ref{Sec: Section 4}–\ref{Sec: Section 10}) contain the detailed proofs of these propositions and related technical results used throughout the argument.

\subsection*{Notation}
We use the symbol $f\lesssim g$ to indicate that $f\leq Cg$ for some constant $C>0$ that may vary from line to line and depends on the main parameters.

\subsection*{Acknowledgements}
The authors thank Lenka Slav\'ikov\'a for suggesting the problem in this paper, her continuous support and helpful comments.

\hfill

\section{Preliminaries}\label{Sec: Section 2}

\subsection{Compactly supported wavelets}

We begin by recalling some standard definitions and results from the theory of compactly supported wavelets.
For any fixed $L\in\bbn$ one can construct real-valued, compactly supported functions $\vartheta_F, \vartheta_M$ in  $\mathscr{C}^L(\bbr)$ satisfying 
  $$\Vert \vartheta_F\Vert_{L^2(\bbr)}=\Vert \vartheta_M\Vert_{L^2(\bbr)}=1, \quad \int_{\bbr}{x^{\alpha}\vartheta_M(x)}dx=0~ \text{for all $0\le \alpha \le L$}$$
For each $k \in \bbz$, $\la\in \bbn_0$, and $G\in \{F,M\}$ we define 
$$
\vartheta_{G,k}^{\la}(x):=2^{\frac{\la}{2}}\vartheta_{G}(2^{\la}x-k).
$$
These functions satisfy the following properties:
\begin{enumerate}
\item[(i)] The family $\{\vartheta_{G,k}^{\lambda}\}_{k\in \mathbb Z}$ has  almost  disjoint  supports. 
\item[(ii)] $\sum_{k\in\bbz }{|\vartheta_{G,k}^{\lambda}(x)|}\lesssim 2^{\frac{\lambda}{2}}$ for all $x \in \mathbb R$.
\end{enumerate}
We specifically assume that the  support of 
 $\vartheta_{G_{j}}$ is contained in $\{\xi\in \bbr: |\xi|\le C_0 \}$ for some $C_0>1$,
which implies 
\begin{equation}\label{supcontgkj}
  \supp(\vartheta_{G_j,k_j}^\la)\subset \big\{\xi_j\in\bbr: |2^{\lambda}\xi_j-k_j|\le C_0 \big\}.
\end{equation}
As a consequence of (i) and (ii) we deduce that for any $0<q<\infty$,
\begin{equation}\label{lrcondition}
  \Big(\sum_{k\in\bbz}\big|\vartheta_{G,k}^{\lambda}(x)\big|^q\Big)^{\frac{1}{q}} \sim_q \sum_{k\in\bbz}\big|\vartheta_{G,k}^{\lambda}(x)\big|\le 2^{\frac{\lambda}{2}}.
\end{equation} 

Now, we define   
$$\Theta_{\GGG}(\xxx):=\vartheta_{G_1}(x_1)\vartheta_{G_{2}}(x_{2})$$
for $\xxx:=(x_1,x_{2})\in \bbr^{2}$ and $\GGG:=(G_1,G_{2})$ in the set     
$$\II:=\big\{\GGG:=(G_1,G_{2}):G_{1},G_2\in\{F,M\} \big\},$$
so that
$$
\Theta_{\GGG,\kkk}^{\la}(\xxx):=2^{\la }\Theta_{\GGG}(2^{\la} \xxx-\kkk)=\vartheta_{G_1,k_1}^{\la}(x_1)\vartheta_{G_2,k_2}^{\la}(x_2).
$$
Then it is known in \cite{Tr2010} that the family of functions
\begin{equation*}
  \bigcup_{\lambda\in\bbn_0}\bigcup_{\kkk\in \bbz^{2}}\big\{ \Theta_{\GGG,\kkk}^{\lambda}(\xxx):\GGG\in \II^{\lambda}\big\}
\end{equation*}
forms an orthonormal basis of $L^2(\bbr^{2})$,
where $\II^0:=\II$ and  for $\lambda \ge 1$, we set  $\II^{\lambda}:= \{(F,G),(G,F),(G,G)\}$.
Moreover, according to \cite[Theorem 1.64]{Tr2006},
if $L$ is sufficiently large, then every  $H\in L^q(\bbr^{2})$ with  $1<q<\infty$ admits the representation
\begin{equation}\label{daubechewavelet}
  H(\xxx)=\sum_{\lambda\in\bbn_0}\sum_{\GGG\in\II^{\lambda}}\sum_{\kkk\in \bbz^{2}}b_{\GGG,\kkk}^{\la}2^{\la }\Theta_{\GGG,\kkk}^{\lambda}( \xxx ),
\end{equation}
where the series converges in $\mathscr S'(\bbr^{2})$, and 
\begin{equation}\label{lqlqs}
  \Big\Vert \Big( \sum_{\la \in \bbn_0} \sum_{\GGG\in\II^{\la}}\sum_{ \vec k\in\bbz^{2}} \big|b_{ \GGG,\kkk}^\la \Theta^\la_{\GGG,\kkk}\big|^2\Big)^{\frac{1}{2}} \Big\Vert_{L^q(\bbr^{2})} \lesssim \|H\|_{L^q(\bbr^{2})}. 
\end{equation}
Here,
\begin{equation*}
  b_{\GGG,\kkk}^{\la}:=\int_{\bbr^{2}}{H(\xxx)\Theta^\la_{ \GGG,\kkk}(\xxx)}d\xxx .
\end{equation*}
By combining (\ref{lqlqs}) with the almost disjoint support property of the functions $\Theta^\la_{ \GGG,\kkk}$, we obtain 
\begin{align}
  \big\Vert \big\{b_{\GGG,\kkk}^{\la}\big\}_{\kkk\in \bbz^{2}}\big\Vert_{\ell^{q}}
  \sim&\Big(2^{\la (2-q)}\int_{\bbr^{2}}\Big( \sum_{\vec k} \big|b^\la_{ \GGG,\kkk}\Theta^\la_{ \GGG,\kkk}(\xxx)\big|^2\Big)^{\frac{q}{2}}d\xxx\Big)^{\frac{1}{q}}\notag \\
  &\lesssim 2^{\la (\frac{2}{q}-1)}\Vert H\Vert_{L^q(\bbr^{2})}.\label{lqestimate}
\end{align}

\subsection{Columns and projections}

For any sets $\UU \subset \bbz^2$, we define the projections
$$\mathcal{P}_1\UU:=\big\{k_1\in\bbz:\kkk=(k_1,k_2) \in\UU ~\text{ for some }k_2\in \bbz \big\},$$
$$\mathcal{P}_2\UU:=\big\{k_2\in\bbz:\kkk=(k_1,k_2) \in\UU ~\text{ for some }k_1\in \bbz \big\}.$$
For each fixed $k_1 \in \mathcal{P}_{1}\UU$,  we define the column set of $\UU$ in the second coordinate by 
\begin{equation*}
  Col_{k_1}^{\UU,2}:=\{k_2 \in \bbz: \kkk=(k_1,k_2) \in\UU\}
\end{equation*}
and, symmetrically, for  $k_2\in\mathcal{P}_2\UU$ we set
\begin{equation*}
  Col_{k_2}^{\UU,1}:=\{k_1 \in \bbz: \kkk=(k_1,k_2) \in\UU\}.
\end{equation*}
With these notations, we can decompose any double sum over $\UU$ as
\begin{equation}\label{disjointdecom1}
  \sum_{\kkk\in \UU}\cdots = \sum_{k_1\in \mathcal{P}_{1}\UU}\Big( \sum_{k_2 \in Col_{k_1}^{\UU,2}}\cdots\Big)=  \sum_{k_2\in \mathcal{P}_{2}\UU}\Big( \sum_{k_1 \in Col_{k_2}^{\UU,1}}\cdots\Big).
\end{equation}
For further discussion of these notations and their higher-dimensional generalizations, we refer to \cite{Gr_He_Ho2018, Gr_He_Ho_Park2023}.

\subsection{Bi-sublinear and weighted interpolation theory}\label{S:mainsectionint}

We begin by recalling a bi-sublinear variant of the Marcinkiewicz interpolation theorem, stated in \cite[Lemma E]{Park2025}, whose proof can be found in \cite[Theorem 1.1]{Gr_Li_Lu_Zh2012}.  

\begin{customlemma}{C}\cite[Lemma E]{Park2025}\label{lem:standint}
  Let $0<p_l^k\leq\infty$ for each $l\in\{1,2\}$, $k=0,1,2,$ and $0<p^{k}\leq\infty$ with $\frac{1}{p^k}=\frac{1}{p_1^k}+\frac{1}{p_2^k}$ for $k=0,1,2$. Suppose that $T$ is an bi-sublinear operator having the mapping properties
  \begin{equation*}
  \|T(f_1,f_2)\|_{L^{p^k,\infty}(\R)}\leq M_k\|f_1\|_{L^{p_1^k}(\R)}\|f_2\|_{L^{p_2^k}(\R)},\qquad k=0,1,2,
  \end{equation*}
  for Schwartz functions $f_1,f_2$ on $\R$. Then for any $0<\theta_k<1$, $k=0,1,2$ with $\theta_0+\theta_1+\theta_2=1$ and $0<p_1,p_2,p\leq\infty$ satisfying
  \begin{equation*}
  \frac{1}{p_1}=\frac{\theta_0}{p_1^0}+\frac{\theta_1}{p_1^1}+\frac{\theta_2}{p_1^2},\qquad \frac{1}{p_2}=\frac{\theta_0}{p_2^0}+\frac{\theta_1}{p_2^1}+\frac{\theta_2}{p_2^2},\qquad\frac{1}{p}=\frac{\theta_0}{p^0}+\frac{\theta_1}{p^1}+\frac{\theta_2}{p^2},
  \end{equation*}
  we have
  \begin{equation*}
   \|T(f_1,f_2)\|_{L^{p,\infty}(\R)}\lesssim M_0^{\theta_0}M_1^{\theta_1}M_2^{\theta_2}\|f_1\|_{L^{p_1}(\R)}\|f_2\|_{L^{p_2}(\R)}.      
  \end{equation*}
  Moreover, if the convex hull of points $(\frac{1}{p_1^0},\frac{1}{p_2^0})$, $(\frac{1}{p_1^1},\frac{1}{p_2^1})$ and $(\frac{1}{p_1^2},\frac{1}{p_2^2})$ forms a non-trivial triangle in $\R^2$, then
  \begin{equation*}
   \|T(f_1,f_2)\|_{L^{p}(\R)}\lesssim M_0^{\theta_0}M_1^{\theta_1}M_2^{\theta_2}\|f_1\|_{L^{p_1}(\R)}\|f_2\|_{L^{p_2}(\R)}.     
   \end{equation*}
\end{customlemma}

The second interpolation result that we need is due to Cao, Olivo and Yabuta \cite[Theorem 3.1]{COY22}. We will use this result in the proof of Theorem \ref{thm:main result} in Subsection \ref{Sec: Section 3.3} below.
\begin{customlemma}{D}\cite[Theorem 3.1]{COY22}\label{lem:COY}
Suppose that $(\Sigma_0, \mu_0),(\Sigma_1, \mu_1)$, $(\Sigma_2, \mu_2)$ and $(\Sigma_3, \mu_3)$ are measure spaces, and $\mathscr{S}_j$ is the collection of all simple functions on $\Sigma_j$, $j=1,2,3$. Denote by $\mathfrak{M}(\Sigma_0)$ the set of all measurable functions on $\Sigma_0$. Let $T: \mathscr{S}=\mathscr{S}_1 \times\mathscr{S}_2\times\mathscr{S}_3 \to \mathfrak{M}(\Sigma_0)$ be a trilinear operator.
Let $0<p_0, q_0< \infty$, $1\le p_j, q_j\le \infty$ $(j=1,2,3)$, and let
$w_j, v_j$ be weights on $\Sigma_j$ $(j=0,1,2,3)$.
Assume that there exist $M_1, M_2 \in (0, \infty)$ such that
\begin{align*}
  &\|T\|_{L^{p_1}(\Sigma_1,\, w_1^{p_1}) \times L^{p_2}(\Sigma_2,\, w_2^{p_2})\times L^{p_3}(\Sigma_3,\, w_3^{p_3}) \to L^{p_0}(\Sigma_0,\, w_0^{p_0})}
  \le M_1,
  \\
  &\|T\|_{L^{q_1}(\Sigma_1,\, v_1^{q_1}) \times
  L^{q_2}(\Sigma_2,\, v_2^{q_2})\times
  L^{q_3}(\Sigma_3,\, v_3^{q_3}) \to L^{q_0}(\Sigma_0,\, v_0^{q_0})} \le
  M_2.
\end{align*}
Then, we have
\begin{equation*}
  \|T\|_{L^{r_1}(\Sigma_1,\, u_1^{r_1}) \times
  L^{r_2}(\Sigma_2,\, u_2^{r_2})\times
  L^{r_3}(\Sigma_3,\, u_3^{r_3}) \to L^{r_0}(\Sigma_0,\, u_0^{r_0})}
  \le M_1^{1-\theta} M_2^{\theta},
\end{equation*}
for all exponents satisfying
\begin{equation*}
  0<\theta<1,\quad \frac{1}{r_j}=\frac{1-\theta}{p_j}+\frac{\theta}{q_j}
  \quad\text{and}\quad u_j=w_j^{1-\theta} v_j^{\theta},\quad j=0,1,2,3.
\end{equation*}
\end{customlemma}

\begin{remark}\label{rmk:COY}
We emphasize that Lemma \ref{lem:COY} also holds for \emph{tri-linearizable} operators, where we recall that a tri-sublinear operator $T$ is called \emph{tri-linearizable} operator if there exist a Banach space $\B$ and a $\B$-valued trilinear operator $U$ such that $T(\vec{f})(x)=\|U(\vec{f})(x)\|_{\B}$. The proof of this fact can be found in \cite[Claim 7.3]{CFRXY24}.
\end{remark}

\hfill

\section{Proof of Theorem \ref{thm:main result}}\label{Sec: Section 3}

The first step of the proof of Theorem \ref{thm:main result} employs the dyadic decomposition of Duoandikoetxea and Rubio de Francia \cite{Du_Ru1986} and utilizes a new reduction step in Subsection \ref{Sec: Section 3.1} that resembles the one as in the proofs of \cite[Theorem 1.2]{Gr_He_Ho_Park2024} and \cite[Theorem 1]{Park2025}. More precisely, let $\Psi$ be a Schwartz function on $\bbr^2$ such that
\begin{equation*}
  \supp{(\wh{\Psi})}\subset \Big\{\xxxi\in\bbr^2: \frac{1}{2}\le |\xxxi|\le 2\Big\}, \q \text{ and }\q \sum_{k\in\bbz}\wh{\Psi_k}(\xxxi)=1,~ \xxxi\not= 0
\end{equation*}
where we set  $\Psi_k:=2^{2k}\Psi(2^k\vec{\cdot}\,)$ for $k\in\bbz$.
For each $\gamma,\mu \in\bbz$,  we define
 $$ K^{\gamma}(\yyy):=\wh{\Psi}(2^{\gamma}\yyy)K(\yyy)~ \text{ and }~ K_{\mu}^{\gamma}(\yyy):=\Psi_{{\mu}+\gamma}\ast K^{\gamma}(\yyy), \quad \yyy\in \bbr^2.$$
 Then $K^\gamma(\yyy)=2^{2\ga } K^0(2^\gamma \yyy)$ and this deduces
 \begin{equation*}
  K_{\mu}^{\gamma}(\yyy)=2^{2\gamma }\big(\Psi_{{\mu}}\ast K^{0}\big)(2^\gamma \yyy)=2^{2\gamma }K_{\mu}^{0}(2^{\gamma}\yyy),
\end{equation*}
or equivalently,
$$
\wh{K^\gamma_\mu}(\xxxi)= \wh{\Psi}(2^{-(\mu+\ga)}\xxxi)\wh{K^0}(2^{-\ga}\xxxi)=\wh{K^0_\mu}(2^{-\ga} \xxxi).
$$
The bilinear associated operator $T_{K_{\mu}^{\gamma}}$ is defined as
$$T_{K_{\mu}^{\gamma}}\big(f_1,f_2\big)(x):=\int_{\bbr^2}K_{\mu}^{\gamma}(y_1,y_2)f_1(x-y_1)f_2(x-y_2)\; d\yyy.$$
so that the bilinear rough singular integral operator $\LL_{\Om}$ can be decomposed as
\begin{equation*}
  \LL_{\Om}=\sum_{\mu\in\Z}\sum_{\gamma\in\Z}T_{K_{\mu}^{\gamma}}.
\end{equation*}

\subsection{Reduction}\label{Sec: Section 3.1}

Let $1<q<\infty$ and we define  $\wh{\Phi}(\yyy):=1-\sum_{\ga\in\bbn}{\wh{\Psi}(2^{-\gamma}\yyy)}$
so that 
\begin{equation*}
  \supp(\wh{\Phi})\subset \{\yyy\in \bbr^2:|\yyy|\le 2\}   
\end{equation*}
and $\wh{\Phi}(\yyy)=1$ for $|\yyy|\le 1$. We first write 
\begin{align}
  \big| \LL_{\Om}^{(2^{\rho})}\big(f_1, f_2\big)(x)\big|
  &\le 
  \bigg| \int_{\bbr^2 }{\Big( K^{(2^{\rho})}(y_1,y_2) - \widetilde{K}^{(2^{\rho})}(y_1,y_2) \Big)f_1(x-y_1)f_2(x-y_2)}\;d\yyy \bigg| \label{ineq 1} \\
  &\qquad \qquad\qquad \qquad+ \bigg| \int_{\bbr^2}   \widetilde{K}^{(2^{\rho})}(y_1,y_2) f_1(x-y_1)f_2(x-y_2) \;d\yyy\bigg|\nonumber,
\end{align}
where
$${K^{(2^\rho)} ( \yyy\,) := K ( \yyy\,)  \chi_{|\yyy\,| > 2^\rho}}, \qq \widetilde{K}^{(2^\rho)} ( \yyy\,) := K ( \yyy\,) \big( 1-  \wh{\Phi} (2^{-\rho}\yyy) \big).$$
Then the first term in the right-hand side of \eqref{ineq 1} is controlled by
$$\mathcal{E}^*\big(f_1,f_2\big)(x):=\sup_{\gamma\in\bbz}\bigg| \int_{\bbr^2}K^{*}_{\gamma}(\yyy)f_1(x-y_1)f_2(x-y_2)\; d\yyy\bigg|,$$
where 
$$K^{*}_{\gamma}(\yyy):= K^{(2^{\gamma})}(\yyy)-\wt{K}^{(2^{\gamma})}(\yyy)=K(\yyy)\wh{\Phi}(2^{-\gamma}\yyy)\chi_{2^{\gamma}<|\yyy|\le 2^{\gamma+1}}(\yyy).$$
We notice that the kernel $K^{*}_{\gamma}$ satisfies
\begin{equation}\label{kbdptest}
  \big|K^{*}_{\gamma}(\yyy)\big|\le \big| K(\yyy)\big|\chi_{2^{\gamma}<|\yyy|\le 2^{\gamma+1}}(\yyy),
\end{equation}
and we have the following estimate.
\begin{lemma}\label{derivativeptest}
  Let $1<q<\infty$ and $\gamma\in\bbz$. For any multi-indices $\alpha$ and for any $0<\delta<\frac{1}{q'}$, we have
\begin{equation*}
  \Big|\partial^{\alpha}\Big( \wh{K^{*}_{\gamma}}(\cdot/2^{\gamma})\Big)(\xxxi)\Big|\lesssim_{q,\delta,\alpha}  \frac{1}{|\xxxi|^{\delta}}\Vert \Omega\Vert_{L^q(\mathbb{S}^1)} 
\end{equation*}
uniformly in $\gamma\in \bbz$.
\end{lemma}
The above lemma is essentially same as \cite[Lemma 2]{Gr_He_Ho2018} and it can be proved by merely mimicking the proof of \cite[Lemma 8.20]{Du2001}.

On the other hand, the second term in the right-hand side of \eqref{ineq 1} is less than
\begin{equation*}
  \bigg|\sum_{\ga\in\bbz: \ga<-\rho}\int_{\bbr^2}K^{\ga}(y_1,y_2) f_1(x-y_1)f_2(x-y_2) \; d\yyy   \bigg|\le \mathcal{L}_{\Omega}^{\sharp}\big(f_1,f_2\big)(x),
\end{equation*}
where
\begin{equation*}
  \LL_{\Om}^{\sharp}\big(f_1,f_2\big)(x) := \sup_{\tau\in \mathbb Z} \Big| \sum_{\gamma<\tau} \sum_{\mu\in\bbz}
  T_{K_{\mu}^{\gamma}}\big(f_1,f_2\big)(x)\Big|.
\end{equation*}

 Therefore, the maximal function $\LL_{\Om}^{\mathrm{dyad},*}(f_1,f_2)$ can be estimated as 
\begin{equation}\label{eq:maindecomp.}
\LL_{\Om}^{\mathrm{dyad},*}\big(f_1,f_2\big)(x)\le \mathcal{E}^*\big(f_1,f_2\big)(x) + \LL_{\Om}^{\sharp}\big(f_1,f_2\big)(x).
\end{equation}

We note that
\begin{equation}\label{emaxomptest}
  \mathcal{E}^*\big(f_1,f_2\big)(x)\lesssim \mathcal{M}_{\Omega}\big(f_1,f_2\big)(x)
\end{equation}
where
\begin{equation*}
  \mathcal M_\Omega\big(f_1, f_2\big)(x)= \sup_{R>0}  \frac{1}{R^{2} } \int_{|\yyy|\le R}   \big|\Omega (\yyy' )    f_1(x-y_1)f_2(x-y_2) \big|  \; d\yyy,
\end{equation*}
  and it is already known in \cite[Lemma 3.1]{Gr_He_Ho_Park2024} that for $1<p,p_1,p_2<\infty$ with $\frac{1}{p}=\frac{1}{p_1}+\frac{1}{p_2}$, and $1<q<\infty$
\begin{equation*}
  \big\Vert \mathcal M_{\Omega}(f_1, f_2)\big\Vert_{L^{p}(\bbr)}\lesssim_q \Vert \Omega\Vert_{L^q(\mathbb S^{1})}
  \Vert f_1\Vert_{L^{p_1}(\bbr)}\Vert f_2\Vert_{L^{p_2}(\bbr)}.
\end{equation*}
This, together with \eqref{emaxomptest}, yields that
\begin{equation*}
  \big\Vert \mathcal{E}^*(f_1, f_2)\big\Vert_{L^{p}(\bbr)}\lesssim_q \Vert \Omega\Vert_{L^q(\mathbb S^{1})}
  \Vert f_1\Vert_{L^{p_1}(\bbr)}\Vert f_2\Vert_{L^{p_2}(\bbr)}.
\end{equation*}
So our main interest is the case when $\frac{1}{2}<p\le 1$.

\hfill

\subsection{Main estimates}\label{Sec: Section 3.2}

We will prove that for all $1<p_1,p_2<\infty$, $\frac{1}{2}<p\le 1$,$1<q<\infty$, and $A>\frac{1}{q}$,
\begin{equation}\label{pmainest1}
  \big\Vert \mathcal{E}^*\big(f_1,f_2\big)\big\Vert_{L^{p}(\bbr)}\lesssim \Vert\Omega\Vert_{L^q(\S^{1},{u_A^q})}\Vert f_1\Vert_{L^{p_1}(\bbr)}\Vert f_2\Vert_{L^{p_2}(\bbr)}
\end{equation}
and 
\begin{equation}\label{pmainest2}
  \big\Vert  \LL_{\Om}^{\sharp}\big(f_1,f_2\big)\big\Vert_{L^{p}(\bbr)}\lesssim \Vert\Omega\Vert_{L^q(\S^{1},{u_A^q})}\Vert f_1\Vert_{L^{p_1}(\bbr)}\Vert f_2\Vert_{L^{p_2}(\bbr)}.
\end{equation}
Then in view of \eqref{eq:maindecomp.}, we conclude
\begin{equation}\label{eq:mainL*Om}
  \big\Vert  \LL_{\Om}^{\mathrm{dyad},*}(f_1,f_2)\big\Vert_{L^{p}(\bbr)}\lesssim \Vert\Omega\Vert_{L^q(\S^{1},u_A^q)}\Vert f_1\Vert_{L^{p_1}(\bbr)}\Vert f_2\Vert_{L^{p_2}(\bbr)}.
\end{equation}

Let us first prove \eqref{pmainest1}. We see
\begin{align}\label{decommathcale}
  \mathcal{E}^*\big(f_1,f_2\big)(x)&\le \sum_{\mu\in\bbz}\sup_{\gamma\in \bbz}\Big|  \LL_{\Omega,\mu}^{\gamma}\big(f_1,f_2\big)(x)  \Big|\nonumber\\
  &=\sum_{\mu\in\bbz: 2^{\mu-10}>C_0} \cdots\qq +\sum_{\mu\in\bbz: 2^{\mu-10}\le C_0}\cdots,
\end{align}
where
\begin{equation*}
  \LL_{\Omega,\mu}^{\gamma}\big(f_1,f_2\big)(x):=\int_{\bbr^2}   \Psi_{\mu-\gamma}\ast K^{*}_{\gamma}(\yyy)f_1(x-y_1)f_2(x-y_2)\; d\yyy   
\end{equation*}
and $C_0$ is the constant that appeared in \eqref{supcontgkj}.
The sums in \eqref{decommathcale} will be treated separately, using the following propositions.

\begin{proposition}\label{muge0tcpg1}
Let $1\le p<\infty$ and $1<p_1,p_2<\infty$ with $\frac{1}{p}=\frac{1}{p_1}+\frac{1}{p_2}$.
Suppose that $\mu\in\bbz$ and $\Omega\in L^1(\S^1)$.
Then we have
$$\Big\Vert \sup_{\gamma\in\bbz}\big| \LL_{\Omega,\mu}^{\gamma}\big(f_1,f_2\big)\big|\Big\Vert_{L^p(\bbr)}\lesssim \Vert \Omega\Vert_{L^1(\mathbb{S}^1)}\Vert f_1\Vert_{L^{p_1}(\bbr)}\Vert f_2\Vert_{L^{p_2}(\bbr)}$$
uniformly in $\mu$.
\end{proposition}

\begin{proposition}\label{mupocasepropo}
Let $1<p_1,p_2<\infty$ and $\frac{1}{2}<p<\infty$ with $\frac{1}{p}=\frac{1}{p_1}+\frac{1}{p_2}$.
Suppose that $\mu\in\Z$ satisfies $2^{\mu-10}\le C_0$ and $\Omega\in L^1(\S^1)$ with $\int_{\mathbb{S}^1}\Omega(\vec{\theta})\; d\sigma(\vec{\theta})=0$.
Then we have
$$\Big\Vert\sup_{\gamma\in \bbz}\big|  \LL_{\Omega,\mu}^{\gamma}\big(f_1,f_2\big)  \big| \Big\Vert_{L^p(\bbr)}\lesssim_{C_0} 2^{2\mu}\Vert \Omega\Vert_{L^1(\mathbb{S}^1)}\Vert f_1\Vert_{L^{p_1}(\bbr)}\Vert f_2\Vert_{L^{p_2}(\bbr)}.$$
\end{proposition}

The above two propositions will be proved in the next two sections.

\begin{proposition}\label{initialpropo}
Let $\frac{4}{3}<q\le 2$ and $\Omega\in L^q(\S^1)$ with $\int_{\mathbb{S}^1}\Omega(\vec{\theta})\; d\sigma(\vec{\theta})=0$. Suppose that $\mu\in\Z$ satisfies $2^{\mu-10}>C_0$.
There exists $\delta>0$, depending on $q$, such that
\begin{equation}\label{initialpropoest}
  \Big\Vert\sup_{\gamma\in\bbz}\big| \LL_{\Omega,\mu}^{\gamma}(f_1,f_2)\big|  \Big\Vert_{L^1(\bbr)}\lesssim_{\delta,C_0} 2^{-\delta\mu}\Vert \Omega\Vert_{L^{q}(\mathbb{S}^1)}\Vert f_1\Vert_{L^2(\bbr)}\Vert f_2\Vert_{L^2(\bbr)}.
\end{equation}
\end{proposition}
This proposition can be proved using the methods employed in \cite{Gr_He_Ho2018, Gr_He_Ho_Park2023}, but since it is quite intricate, we will provide it in a tailored setting in Section \ref{Sec: Section 6} below.

\begin{proposition}\label{PropoHLS}
Suppose that $\mu\in\Z$ satisfies $2^{\mu-10}>C_0$ and $\Omega\in L^1(\mathbb{S}^1,u_1)$, where we recall $u_1(\theta_1,\theta_2)=\frac{1}{|\theta_1-\theta_2|}$ for $(\theta_1,\theta_2)\in \mathbb{S}^1$. Assume that $\int_{\mathbb{S}^1}\Omega(\vec{\theta})\; d\sigma(\vec{\theta})=0$.
For any $0<\epsilon<1$,
$$\bigg\Vert \sup_{\gamma\in\bbz}\big| \LL_{\Omega,\mu}^{\gamma}(f_1,f_2)\big|\bigg\Vert_{L^{\frac{1}{2},\infty}(\R)}\lesssim_{\epsilon,C_0}2^{\epsilon \mu}\Vert \Omega\Vert_{L^1(\mathbb{S}^1,{u_1})}\Vert f_1\Vert_{L^1(\bbr)}\Vert f_2\Vert_{L^1(\bbr)}.$$
\end{proposition}
The proof is based on the idea in \cite[Proof of Proposition 4]{Ho_La_Sl_submitted} and it will be given in Section \ref{Sec: Section 7}.

By recalling the embeddings $L^p(\mathbb R) \hookrightarrow L^{p,\infty}(\mathbb R)$ and $L^\infty(\mathbb \S^1) \hookrightarrow L^{q}(\mathbb \S^1)$ for all $1\leq q<\infty$ and interpolating between the estimates in Propositions \ref{muge0tcpg1} and \ref{initialpropo} via  Lemma \ref{lem:standint}, we obtain that for any $1< p<\infty$ and $1<p_1,p_2<\infty$, there exists $\delta_1>0$ such that
\begin{equation}\label{generalrgestimp}
\Big\|  \sup_{\gamma\in\bbz}\big| \LL_{\Omega,\mu}^{\gamma}\big(f_1,f_2\big)\big| \Big\|_{L^{p}(\bbr) } \lesssim 2^{-\delta_1 \mu}\Vert \Omega\Vert_{L^{\infty}(\mathbb{S}^{1})}\Vert f_1\Vert_{L^{p_1}(\bbr)}\Vert f_2\Vert_{L^{p_2}(\bbr)},
\end{equation}
when $\mu\in\bbz$ with $2^{\mu-10}>C_0$.

By applying the same argument in \cite[Section 3]{Park2025}, we obtain the following improved estimates, where the $L^{\infty}$ norm of $\Omega$ in \eqref{generalrgestimp} by the $L^q$ norm for any $q>1$.
\begin{proposition}\label{prop:intmaxE}
Suppose that $\mu\in\Z$ satisfies $2^{\mu-10}>C_0$ and let $1<p_1,p_2<\infty$, $1< p<\infty$ satisfy $\frac{1}{p}=\frac{1}{p_1}+\frac{1}{p_2}$. Assume that $q>1$, $\Omega\in L^q(\S^1)$ with $\int_{\mathbb{S}^1}\Omega(\vec{\theta})\, d\sigma(\vec{\theta})=0$. Then there exists $\delta>0$ depending on $p_1,p_2$ and $q$ such that
\begin{equation*}
  \Big\|\sup_{\gamma\in\bbz}\big|\LL_{\Omega,\mu}^{\gamma}\big(f_1,f_2\big)\big| \Big\|_{L^{p}(\bbr) } \lesssim 2^{-\delta \mu}\Vert \Omega\Vert_{L^{q}(\mathbb{S}^{1})}\Vert f_1\Vert_{L^{p_1}(\bbr)}\Vert f_2\Vert_{L^{p_2}(\bbr)}.    
\end{equation*}
\end{proposition}
The proposition will be proved in Section \ref{Sec: Section 8}.

Now, with the aid of Lemma \ref{lem:standint} we interpolate further between the estimates in Propositions \ref{PropoHLS} and \ref{prop:intmaxE}. In particular, we obtain the following exponential decay estimate in $\mu$ related to the family of maximal operators  $\sup_{\gamma\in\bbz}\big| \LL_{\Omega,\mu}^{\gamma}|$.

\begin{proposition}\label{prop:intmaxE1}
Suppose that $\mu\in \Z$ satisfies $2^{\mu-10}>C_0$ and let $1<p_1,p_2<\infty$, $\frac{1}{2}<p \le 1$ satisfy $\frac{1}{p}=\frac{1}{p_1}+\frac{1}{p_2}$. Assume that  $q>1$ and  $\Omega \in L^q(\S^{1},u_A^q)$ for $A>\frac{1}{q}$ with $\int_{\S^{1}} \Omega(\vec\theta) \,d\sigma(\vec\theta)=0$. Then there exists $\delta>0$ such that
\begin{equation*}
  \Big\|\sup_{\gamma\in\bbz}\big| \LL_{\Omega,\mu}^{\gamma}\big(f_1,f_2\big)\big| \Big\|_{L^{p}(\bbr) }\lesssim_{p_1,p_2,q,A}
  2^{-\delta \mu}\Vert\Omega\Vert_{L^q(\S^{1},{u_A^q})}\Vert f_1\Vert_{L^{p_1}(\R)}\Vert f_2\Vert_{L^{p_2}(\R)}.
\end{equation*}   
\end{proposition}

The proposition will be proved in Section \ref{Sec: Section 9}.

Now, let us complete the proof of \eqref{pmainest1}.
Recalling \eqref{decommathcale},
\begin{align}\label{eq:decommathcale1}
  \big\Vert\mathcal{E}^*(f_1,f_2)\big\Vert_{L^p(\R)}&\leq\bigg(\sum_{\mu\in\bbz:2^{\mu-10}>C_0}\Big\Vert \sup_{\gamma\in\bbz}\big| \LL_{\Omega,\mu}^{\gamma}\big(f_1,f_2\big)\big|\Big\Vert_{L^p(\bbr)}^p\\
  &\qquad+\sum_{\mu\in\bbz:2^{\mu-10}\leq C_0}\Big\Vert \sup_{\gamma\in\bbz}\big| \LL_{\Omega,\mu}^{\gamma}\big(f_1,f_2\big)\big|\Big\Vert_{L^p(\bbr)}^p\bigg)^{\frac{1}{p}}\notag.
\end{align}
Using Proposition \ref{mupocasepropo} we know that for $\mu\in\Z$ which satisfies $2^{\mu-10}\leq C_0$ we have
\begin{equation*}
  \Big\Vert \sup_{\gamma\in\bbz}\big| \LL_{\Omega,\mu}^{\gamma}\big(f_1,f_2\big)\big|\Big\Vert_{L^p(\bbr)} \lesssim 2^{2\mu}\Vert \Omega\Vert_{L^1(\mathbb{S}^1)}\Vert f_1\Vert_{L^{p_1}(\bbr)}\Vert f_2\Vert_{L^{p_2}(\bbr)},
\end{equation*}
up to a constant depending on $p_1$, $p_2$, $q$ and $C_0$. Hence,
\begin{equation}\label{eq:decommathcale2}
  \sum_{\mu\in\bbz:2^{\mu-10}\leq C_0}\Big\Vert \sup_{\gamma\in\bbz}\big| \LL_{\Omega,\mu}^{\gamma}\big(f_1,f_2\big)\big|\Big\Vert_{L^p(\bbr)}\lesssim\Vert\Omega\Vert_{L^q(\S^{1},u_A^q)}\Vert f_1\Vert_{L^{p_1}(\bbr)}\Vert f_2\Vert_{L^{p_2}(\bbr)},
\end{equation}
for $A>\frac{1}{q}$, where we applied the embedding $L^q(\S^{1},u_A^q)\hookrightarrow L^{1}(\S^{1})$.

For $\mu\in\Z$ which satisfies $2^{\mu-10}>C_0$,
by recalling Proposition \ref{prop:intmaxE1}, we know that 
there exists a small constant $\delta_2=\delta_2(p_1,p_2,q)>0$, such that
\begin{equation*}
  \Big\Vert\sup_{\gamma\in\bbz}\big| \LL_{\Omega,\mu}^{\gamma}\big(f_1,f_2\big)\big|\Big\Vert_{L^p(\bbr)} \lesssim2^{-\mu\delta_2} \Vert\Omega\Vert_{L^q(\S^{1},u_A^q)}\Vert f_1\Vert_{L^{p_1}(\R)}\Vert f_2\Vert_{L^{p_2}(\R)},
\end{equation*}
which gives that 
\begin{equation}\label{eq:decommathcale3}
  \sum_{\mu\in\bbz:2^{\mu-10}>C_0}\Big\Vert\sup_{\gamma\in\bbz}\big| \LL_{\Omega,\mu}^{\gamma}\big(f_1,f_2\big)\big|\Big\Vert_{L^p(\bbr)}\lesssim \Vert\Omega\Vert_{L^{q}(\S^1,{u_A^q})}\Vert f_1\Vert_{L^{p_1}(\R)}\Vert f_2\Vert_{L^{p_2}(\R)}.    
\end{equation}
Hence, by combining  \eqref{eq:decommathcale1}, \eqref{eq:decommathcale2} and \eqref{eq:decommathcale3} we conclude the desired estimate \eqref{pmainest1}. \\

\hfill

In order to prove \eqref{pmainest2},
we write
\begin{align}\label{eq:maxsharp}
  \LL_{\Omega}^{\sharp}\big(f_1,f_2\big)(x)&\le \sup_{\tau\in\bbz}\Big| \sum_{\gamma<\tau}\sum_{\mu\in\bbz: 2^{\mu-10}\le C_0}T_{K_{\mu}^{\gamma}}\big(f_1,f_2\big)(x)\Big|\\
  &\qquad+\sup_{\tau\in\bbz}\Big| \sum_{\gamma<\tau}\sum_{\mu\in\bbz:2^{\mu-10}>C_0}T_{K_{\mu}^{\gamma}}\big(f_1,f_2\big)(x)\Big|\notag    
\end{align}
and it is known in \cite[Proposition 4.1]{Gr_He_Ho_Park2024} that for any $1<p_1,p_2<\infty$ with $\frac{1}{p}=\frac{1}{p_1}+\frac{1}{p_2}$ and for any $1<q<\infty$
\begin{equation}\label{eq:max0}
  \bigg\Vert \sup_{\tau\in\bbz}\Big| \sum_{\gamma<\tau}\sum_{\mu\in\bbz: 2^{\mu-10}\le C_0}T_{K_{\mu}^{\gamma}}\big(f_1,f_2\big)\Big|  \bigg\Vert_{L^p(\bbr)} \lesssim \Vert \Omega\Vert_{L^q(\mathbb{S}^1)}\Vert f_1\Vert_{L^{p_1}(\bbr)}\Vert f_2\Vert_{L^{p_2}(\bbr)}. 
\end{equation}

To deal with the other term, we define
\begin{equation}\label{eq:sharpomegamu}
  \LL_{\Omega,\mu}^{\sharp}\big(f_1,f_2\big)(x):=\sup_{\tau\in\bbz}\bigg| \sum_{\gamma<\tau}T_{K_{\mu}^{\gamma}}\big(f_1,f_2\big)(x)\bigg|    
\end{equation}
for each $\mu\in\bbn$.
Then
\begin{equation*}
  \sup_{\tau\in\bbz}\Big| \sum_{\gamma<\tau}\sum_{\mu\in\bbz:2^{\mu-10}>C_0}T_{K_{\mu}^{\gamma}}\big(f_1,f_2\big)(x)\Big| \le \sum_{\mu\in\bbz:2^{\mu-10}>C_0}\LL_{\Omega,\mu}^{\sharp}\big(f_1,f_2\big)(x).    
\end{equation*}
We claim that for $1<p_1,p_2<\infty$, $\frac{1}{2}<p<\infty$, and $q>1$, there exists $\delta>0$
\begin{equation}\label{eq:int1}
  \big\Vert \LL_{\Omega,\mu}^{\sharp}(f_1,f_2)\big\Vert_{L^{p}(\R)}\lesssim
  2^{-\delta\mu}\Vert\Omega\Vert_{L^q(\S^{1},u_A^q)}\Vert f_1\Vert_{L^{p_1}(\R)}\Vert f_2\Vert_{L^{p_2}(\R)}.
\end{equation}
Then in view of \eqref{eq:maxsharp} and \eqref{eq:max0}, the desired estimate \eqref{pmainest2} follows.\\

Now, let us prove \eqref{eq:int1}.
It is known in \cite[(1.12)]{Park2025} that for any $1<p_1,p_2<\infty$ with $\frac{1}{p}=\frac{1}{p_1}+\frac{1}{p_2}$, there exists a $\delta>0$ such that
\begin{equation}\label{prop:sharpmax1}
  \big\Vert \LL_{\Omega,\mu}^{\sharp}\big(f_1,f_2\big) \big\Vert_{L^{p}(\R)}\lesssim_{\delta}2^{-\delta \mu}\Vert \Omega\Vert_{L^q(\S^1)}\Vert f_1\Vert_{L^{p_1}(\bbr)}\Vert f_2\Vert_{L^{p_2}(\bbr)},
\end{equation}
provided that 
\begin{equation}\label{qcondi}
  q>\max\Big\{1,\frac{1}{2-\frac{1}{p}}\Big\}.
\end{equation}
So when $p\ge 1$, we are done as \eqref{qcondi} becomes $q>1$. For the remaining case $\frac{1}{2}<p<1$, we apply an interpolation method with the following end-point estimates.

\begin{proposition}\label{lommushweakre}
Assume that $\Omega\in L^1(\mathbb{S}^1,u_1)$ with $\int_{\mathbb{S}^1}\Omega(\vec{\theta})\; d\sigma(\vec{\theta})=0$. For any $\epsilon>0$, we have
$$\big\Vert \LL_{\Omega,\mu}^{\sharp}\big(f_1,f_2\big) \big\Vert_{L^{\frac{1}{2},\infty}(\bbr)}\lesssim_{\epsilon}2^{\epsilon\mu}\Vert \Omega\Vert_{L^1(\mathbb{S}^1,{u_1})}\Vert f_1\Vert_{L^{1}(\bbr)}\Vert f_2\Vert_{L^{1}(\bbr)}.$$
\end{proposition}
The proof of the proposition will be given in Section \ref{Sec: Section 10}. We notice that the inequality in Proposition \ref{lommushweakre} itself is not enough good to estimate the $L^p$ norm of the left-hand side of \eqref{eq:sharpomegamu} because the upper bound is not summable for $\mu\in \bbz$ with $2^{\mu-10}>C_0$ due to the exponential growth $2^{\epsilon \mu}$. However, an important observation is that $\epsilon$ can be chosen arbitrarily small, which in turn allows the exponential growth rate to be made as small as desired. \\

Now, we fix $1<p_1,p_2<\infty$ and $\frac{1}{2}<p< 1$ satisfying $\frac{1}{p}=\frac{1}{p_1}+\frac{1}{p_2}$. By choosing $\eta>0$ small enough depending on $p_1$, $p_2$, we may conclude that the point $A_1:=(\frac{1}{p_1},\frac{1}{p_2})$ lies in the interior of a non-trivial triangle in $\R^2$ which is formed by the convex hull of points $A_2:=(1,1)$, $A_3:=(1-2\eta,\eta)$ and $A_4:=(\eta, 1-2\eta)$. See the following Figure \ref{fig:figure2} for the interpolation.

%%%%%%%%%%%%%%%%%%%%

\begin{figure}[h]
\begin{tikzpicture}
\path[fill=red!05] (0.1,3.8)--(3.8,0.1)--(4,4)--(0.1,3.8);
\draw [<->] (0,5)--(0,0)--(5,0);
\draw[dash pattern= { on 2pt off 1pt}] (0,4)--(4,4)--(4,0);
\draw[-] (0,4)--(4,0);%
\node [below left] at (0,0) {\tiny$(0,0)$};
\filldraw[fill=black] (3.8,0.1)  circle[radius=0.3mm];%%A_3
\filldraw[fill=black] (2.7,2.7)  circle[radius=0.5mm];%%A_1
\filldraw[fill=black] (0.1,3.8)  circle[radius=0.3mm];%%A_4
\filldraw[fill=black] (4,4)  circle[radius=0.3mm];%%A_2
\draw [red, dash pattern= { on 2pt off 1pt}] (0.1,3.8)--(3.8,0.1)--(4,4)--(0.1,3.8);%%A_4A_3A_2A_4
\draw[-] (2.7,2.7)--(4.5,3);
\node [right] at (4.5,3) {\tiny$A_1=(\frac{1}{{p_1}},\frac{1}{{p_2}})$};%
\draw[-] (3.8,0.1)--(4.5,0.6);%%
\node [right] at (4.5,0.6) {\tiny$A_3=(1-2\eta,\eta)$};%%
\draw[-] (0.9,4.4)--(0.1,3.8);%
\node [above] at (1.4,4.4) {\tiny$A_4=(\eta,1-2\eta) $};%%
\node [below] at (4,0) {\tiny$(1,0)$};
\node [left] at (0,4) {\tiny$(0,1)$};
\node [above right] at (4,4) {\tiny$A_2=(1,1) $};%%
\node [below] at (5,0) {\tiny$\frac{1}{p_1}$};
\node [left] at (0,5) {\tiny${\frac{1}{p_2}}$};
\end{tikzpicture}
\caption{Interpolation between estimates at points $A_2$, $A_3$ and $A_4$}\label{fig:figure2}
\end{figure}

%%%%%%%%%%%%%%%%%%%%

Hence, there exist  $0<\theta_0,\theta_1,\theta_2<1$ with 
\begin{equation*}
  \theta_0+\theta_1+\theta_2=1
\end{equation*}
such that 
\begin{align}
  \frac{1}{p_1}&=\theta_0+(1-2\eta)\theta_1+\eta\theta_2,\label{eq:main eq 3}\\
  \frac{1}{p_2}&=\theta_0+\eta\theta_1+(1-2\eta)\theta_2\label{eq:main eq 4}.
\end{align}
By adding \eqref{eq:main eq 3} and~\eqref{eq:main eq 4}, we get
\begin{equation}\label{eq:main eq 5}
  \frac{1}{p}=2\theta_0+(1-\eta)\theta_1+(1-\eta)\theta_2.
\end{equation}

Let $q>1$. As in the proofs of Proposition \ref{prop:intmaxE1} or \cite[Proposition 6]{Ho_La_Sl_submitted} we recall the embeddings $L^q(\S^{1},u_A^q)\hookrightarrow L^{q}(\S^{1})$ and $L^q(\S^{1},u_A^q)\hookrightarrow L^1(\S^{1},u_1)$ for $A>\frac{1}{q}$.
Having at our disposal these embeddings, it follows from Proposition \ref{lommushweakre} and the estimate \eqref{prop:sharpmax1} that for arbitrary $\epsilon>0$ and for some $\wt{\delta}>0$,
\begin{align*}
  \big\Vert \LL_{\Omega,\mu}^{\sharp}\big\Vert_{L^1(\R)\times L^1(\R)\to L^{\frac{1}{2},\infty}(\R)}&\lesssim_{\epsilon} 2^{\epsilon\mu}\Vert\Omega\Vert_{L^q(\S^{1},u_A^q)},\\
  \big\Vert \LL_{\Omega,\mu}^{\sharp}\big\Vert_{L^\frac{1}{1-2\eta}(\R)\times L^\frac{1}{\eta}(\R)\to L^{\frac{1}{1-\eta},\infty}(\R)}&\lesssim 2^{-{\wt{\delta}}\mu}\Vert\Omega\Vert_{L^q(\S^{1},{u_A^q})},\\
  \big\Vert \LL_{\Omega,\mu}^{\sharp}\big\Vert_{L^\frac{1}{\eta}(\R)\times L^\frac{1}{1-2\eta}(\R)\to L^{\frac{1}{1-\eta},\infty}(\R)}&\lesssim 2^{-{\wt{\delta}}\mu}\Vert\Omega\Vert_{L^q(\S^{1},u_A^q)},
\end{align*}
up to constants depending on $\eta$, $\epsilon$, $q$, $\wt{\delta}$, and $A$.
Then applying Lemma \ref{lem:standint} we deduce that
\begin{equation}\label{lsharpommuest}
  \big\Vert \LL_{\Omega,\mu}^{\sharp}\big\Vert_{L^{p_1}(\R)\times L^{p_2}(\R)\to L^{p}(\R)}\lesssim 
  2^{-\mu(\wt{\delta}\theta_1+\wt{\delta}\theta_2-\epsilon \theta_0)     }\Vert\Omega\Vert_{L^q(\S^{1},u_A^q)},
\end{equation}
up to a constant depending on $p_1$, $p_2$, $q$, $\epsilon$ and $A$.
Choosing $0<\epsilon<\wt{\delta} \frac{\theta_1+\theta_2}{\theta_0}$ so that $ \wt{\delta} \theta_1+\wt{\delta} \theta_2-\epsilon\theta_0>0$ and setting $\delta:= \wt{\delta} \theta_1+\wt{\delta} \theta_2 -\epsilon \theta_0$, we finally derive inequality~\eqref{eq:int1}.

\hfill

\subsection{Improved estimates via interpolation}\label{Sec: Section 3.3} 

Now, let us  finish the proof of Theorem \ref{thm:main result} by improving the estimate \eqref{eq:mainL*Om} via interpolation.

Let $1<p_1,p_2<\infty$ and $\frac{1}{2}<p\le 1$ with $\frac{1}{p}=\frac{1}{p_1}+\frac{1}{p_2}$, and let $q>1$ satisfy $\frac{1}{p}+\frac{1}{q}\geq 2$. 
Assuming \eqref{sharpcond}, we choose $0<\eta<\frac{1}{2}$ sufficiently small so that 
\begin{equation*}
  A>\frac{1}{p}+\frac{1}{q}-2+2\eta+\frac{\eta}{1+\eta}=:A_0>\frac{1}{p}+\frac{1}{q}-2.
\end{equation*}
Now, it is possible to select $\wt{q}>1$ and $\eta<\kappa<1-\eta$ so that the point $(\frac{1}{q}, \frac{1}{p_1},\frac{1}{p_2})$ lies on the line segment in $\R^3$ connecting the two points $(\frac{1}{\widetilde{q}},1-\eta, 1-\eta)$ and $(\frac{1}{1+\eta}, \kappa-\eta, 1-\kappa-\eta)$.

In particular, setting
\begin{align*}
  \widetilde{q}&=\frac{(1+\eta)(-1+2\eta+\frac{1}{p})}{-2+\frac{1}{p}+\frac{1}{q}+2\eta+\frac{\eta}{q}},\\
  \kappa&=\frac{1}{2-2\eta-\frac{1}{p}}\Big(\eta-1+\frac{1}{p_1}\Big)+1,\\
  \theta&=2-2\eta-\frac{1}{p},
\end{align*}
the following equations are satisfied:
\begin{align*}
  \frac{1}{q}&=\frac{1-\theta}{\wt{q}}+\frac{\theta}{1+\eta}\\
  \frac{1}{p_1}&=(1-\eta)(1-\theta)+(\kappa-\eta)\theta,\\
  \frac{1}{p_2}&=(1-\eta)(1-\theta)+(1-\kappa-\eta)\theta,\\
  \frac{1}{p}&=(2-2\eta)(1-\theta)+(1-2\eta)\theta.
\end{align*}

Before applying interpolation method, we notice that the maximal operator $\LL_{\Om}^{\mathrm{dyad},*}$ is bi-sublinear and therefore we cannot yet use the weighted interpolation in Lemma \ref{lem:COY}. To overcome this obstacle, we view $\LL_{\Om}^{\mathrm{dyad},*}$ as a tri-linearizable operator (see detailed definition in Remark \ref{rmk:COY} of Subsection \ref{S:mainsectionint}). To be more precise, for each $\rho\in\bbz$, we define
\begin{equation*}
  \TT^{\rho}\big(\Omega,f_1,f_2\big)(x):=\int_{|\yyy|>2^{\rho}} \frac{\left(\Omega - 2Avg(\Omega) \chi_{\{(\theta_1,\theta_2): \theta_1 \theta_2 <0\}}\right)(\yyy')}{|\yyy|^{2}}\prod_{j=1}^2 f_j(x-y_j)\,d\yyy,
\end{equation*}
where $Avg(\Omega) =\frac{1}{2\pi}\int_{\S^{1}} \Omega \,d\s$ is the average of $\Omega$,
so that
$\big\{ \TT^{\rho}\big\}_{\rho\in\bbz}$ is a $\ell^{\infty}(\bbz)$-valued trilinear operator and let
\begin{equation*}
  \TT^{*}(\Omega, f_1,f_2)(x):=\big\Vert \big\{ \TT^{\rho}\big(\Omega,f_1,f_2\big)(x)\big\}_{\rho\in\bbz} \big\Vert_{\ell^{\infty}(\bbz)}.
\end{equation*}
Then $\TT^*$ is a tri-linearizable operator and we notice that $\TT^*(\Omega, f_1,f_2)$ agrees with $\LL_{\Om}^{\mathrm{dyad}*}(f_1,f_2)$ for every $\Omega$ with mean value zero. Thus, according to Remark \ref{rmk:COY} we are now able to proceed with our interpolation argument.    

Observing $1+\eta, \frac{1}{\kappa-\eta}, \frac{1}{1-\kappa-\eta},\frac{1}{1-\eta}>1$, and setting 
$$\widetilde{A}=\frac{(1+\eta)(-2+\frac{1}{p}+\frac{1}{q}+2\eta)+\eta}{(1+\eta)(-1+2\eta+\frac{1}{p})}\bigg(>\frac{-2+\frac{1}{p}+\frac{1}{q}+2\eta+\frac{\eta}{q}}{(1+\eta)(-1+2\eta+\frac{1}{p})}=\frac{1}{\widetilde{q}}\bigg)$$
so that $$\wt{A}(1-\theta)=A_0,$$
it follows from Theorem \ref{thm:Park2025} and \eqref{eq:mainL*Om} that

\begin{align*}
  &\big\Vert \TT^*\big\Vert_{L^{1+\eta}(\S^{1}) \times L^\frac{1}{\kappa-\eta}(\R)\times L^\frac{1}{1-\kappa-\eta}(\R)\to L^{\frac{1}{1-2\eta}}(\R)}\lesssim 1,\\
  &\big\Vert \TT^*\big\Vert_{L^{\tilde{q}}(\S^{1},u_{\wt{A}}^{\tilde{q}}) \times L^\frac{1}{1-\eta}(\R)\times L^\frac{1}{1-\eta}(\R)\to L^{\frac{1}{2-2\eta}}(\R)}\lesssim 1.
\end{align*}
Now, based on Remark \ref{rmk:COY}, we apply Lemma \ref{lem:COY} in order to get
\begin{equation*}
  \big\Vert \TT^*\big\Vert_{L^q(\S^{1},u_{A_0}^q) \times L^{p_1}(\R)\times L^{p_2}(\R)\to L^{p}(\R)}\lesssim1.
\end{equation*}
Since $u_{A_0}(\theta_1,\theta_2)\lesssim u_A(\theta_1,\theta_2)$, the desired result finally follows. This completes the proof of Theorem \ref{thm:main result}.

\hfill

\section{Proof of Proposition \ref{muge0tcpg1}}\label{Sec: Section 4}

Using \eqref{kbdptest},
\begin{align*}
 \big| \LL_{\Omega,\mu}^{\gamma}\big(f_1,f_2\big)(x)\big|&\le \int_{\bbr^2}   \big| \Psi_{\mu-\gamma}\big|\ast |K^{*}_{\gamma}|(\yyy)\big| f_1(x-y_1)f_2(x-y_2)\big|\; d\yyy\\
 &\le \int_{|\zzz|\sim 2^{\gamma}}\big| K(\zzz)\big| \Big(\int_{\bbr^2}\big| \Psi_{\mu-\gamma}(\yyy-\zzz)f_1(x-y_1)f_2(x-y_2)\big|\; d\yyy\Big)\; d\zzz\\
 &\lesssim \int_{|\zzz|\sim 2^{\gamma}}\big| K(\zzz)\big|\mathcal{M}f_1(x-z_1)\mathcal{M}f_2(x-z_2)\; d\zzz\\
 &\lesssim \int_{\mathbb{S}^1}\big| \Omega(\vec{\theta})\big|\;\mathbb{M}^{\theta_1}\mathcal{M}f_1(x)\;\mathbb{M}^{\theta_2}\mathcal{M}f_2(x)\; d\sigma(\vec{\theta})
\end{align*}
where the directional maximal operator $\mathbb{M}^{\theta}$ is defined by
$$\mathbb{M}^{\theta_j}g(x):=\sup_{r>0}\Big( \frac{1}{r}\int_0^r\big| g(x-t\theta_j)\big|\;dt\Big).$$
Therefore
\begin{align*}
  \Big\Vert\sup_{\gamma\in\bbz}\big| \LL_{\Omega,\mu}^{\gamma}\big(f_1,f_2\big)\big|\Big\Vert_{L^p(\bbr)}&\lesssim \int_{\mathbb{S}^1} \big| \Omega(\vec{\theta})\big|  \; \big\Vert \mathbb{M}^{\theta_1}\mathcal{M}f_1 \;\mathbb{M}^{\theta_2}\mathcal{M}f_2\big\Vert_{L^p(\bbr)}  \; d\sigma(\vec{\theta})\\
  &\le \int_{\mathbb{S}^1} \big| \Omega(\vec{\theta})\big|  \; \big\Vert \mathbb{M}^{\theta_1}\mathcal{M}f_1 \big\Vert_{L^{p_1}(\bbr)}\big\Vert \mathbb{M}^{\theta_2}\mathcal{M}f_2\big\Vert_{L^{p_2}(\bbr)}  \; d\sigma(\vec{\theta})\\
  &\lesssim \Vert \Omega\Vert_{L^1(\mathbb{S}^1)}\Vert f_1\Vert_{L^{p_1}(\bbr)}\Vert f_2\Vert_{L^{p_2}(\bbr)},
\end{align*}
where the last estimate follows from the $L^{p_j}$ boundedness of $\mathbb{M}^{\theta_j}$ and $\mathcal{M}$ with constants independent of $\theta_j$ (by the method of rotation).
This completes the proof.

\hfill

\section{Proof of Proposition \ref{mupocasepropo}}\label{Sec: Section 5}

We will actually prove that for $\mu\in\bbz$ with $2^{\mu-10}\le C_0$,
\begin{equation}\label{minorptmest}
 \big| \LL_{\Omega,\mu}^{\gamma}\big(f_1,f_2\big)(x)\big|\lesssim 2^{2\mu}\Vert \Omega\Vert_{L^1(\mathbb{S}^1)}\mathcal{M}f_1(x)\mathcal{M}f_2(x) \q \text{uniformly in }~\gamma.
\end{equation}
Then the desired result follows immediately from the $L^{q}$ boundedness of $\mathcal{M}$ for $q>1$.

To verify \eqref{minorptmest}, we write
\begin{align*}
 \big| \LL_{\Omega,\mu}^{\gamma}\big(f_1,f_2\big)(x)\big|&\le \int_{|\yyy|\le 2^{\gamma+10}}   \big| \Psi_{\mu-\gamma}\ast K^{*}_{\gamma}(\yyy)\big|\big| f_1(x-y_1)f_2(x-y_2)\big|\; d\yyy\\
 &\qq+ \int_{|\yyy|>2^{\gamma+10}}   \big| \Psi_{\mu-\gamma}\ast K^{*}_{\gamma}(\yyy)\big| f_1(x-y_1)f_2(x-y_2)\big|\; d\yyy\\
 &=:\mathfrak{J}_{\mathrm{in}}+\mathfrak{J}_{\mathrm{out}}
\end{align*}
We first apply \eqref{kbdptest} to obtain
\begin{align*}
  \big| \Psi_{\mu-\gamma}\ast K^{*}_{\gamma}(\yyy)\big|\le 2^{2\mu-2\gamma}\int_{|\zzz|\sim 2^{\gamma}}\big| K(\zzz)\big|\; d\zzz \lesssim 2^{2\mu-2\gamma}\Vert \Omega\Vert_{L^1(\mathbb{S}^1)},
\end{align*}
Then
\begin{align*}
  \mathfrak{J}_{\mathrm{in}}&\lesssim 2^{2\mu}\Vert \Omega\Vert_{L^1(\mathbb{S}^1)} \Big( 2^{-2\gamma}\int_{|\yyy|\le 2^{\gamma+10}} \big|f_1(x-y_1)f_2(x-y_2) \big|\;d\yyy    \Big)\\
  &\lesssim 2^{2\mu}\Vert \Omega\Vert_{L^1(\mathbb{S}^1)} \mathcal{M}f_1(x)\mathcal{M}f_2(x).
\end{align*}
as expected.
Moreover, using  \eqref{vanishingmtcondition} and \eqref{kbdptest},
\begin{align*}
  \big|\Psi_{\mu-\gamma}\ast K^{*}_{\gamma}(\yyy)\big|&\le \int_{\bbr^2}\big| K^{*}_{\gamma}(\zzz)\big|\Big| \Psi_{\mu-\gamma}(\yyy-\zzz)-\Psi_{\mu-\gamma}(\yyy)\Big|\; d\zzz\\
  &\le\int_{2^{\gamma}<|\zzz|\le 2^{\gamma+1}}\big| K(\zzz)\big| \Big( 2^{\mu-\gamma}|\zzz| \int_0^1 \frac{2^{2(\mu-\gamma)}}{(1+2^{\mu-\gamma}|\yyy-s\zzz|)^L}\; ds\Big) \; d\zzz
\end{align*}
for $L>2$.
If $|\yyy|>2^{\gamma+10}$, $2^{\gamma}<|\zzz|\le 2^{\gamma+1}$, and $0<s<1$, then
$$|\yyy-s\zzz|\gtrsim |\yyy|$$
and thus we obtain
\begin{equation*}
  \big| \Psi_{\mu-\gamma}\ast K^{*}_{\gamma}(\yyy)\big|\lesssim  2^{\mu} \frac{2^{2(\mu-\gamma)}}{(1+2^{\mu-\gamma}|\yyy|)^L}\int_{|\zzz|\sim 2^{\gamma}}|K(\zzz)|\; d\zzz\lesssim  2^{\mu} \Vert \Omega\Vert_{L^1(\mathbb{S}^1)}\frac{2^{2(\mu-\gamma)}}{(1+2^{\mu-\gamma}|\yyy|)^L}.
\end{equation*}
This proves
\begin{align*}
  \mathfrak{J}_{\mathrm{out}}&\lesssim 2^{\mu} \Vert \Omega\Vert_{L^1(\mathbb{S}^1)}\int_{\bbr^2}   \frac{2^{2(\mu-\gamma)}}{(1+2^{\mu-\gamma}|\yyy|)^L}\big| f_1(x-y_1) f_2(x-y_2)\big|    \;d\yyy\\
  &\lesssim 2^{\mu}\Vert \Omega\Vert_{L^1(\mathbb{S}^1)}\mathcal{M}f_1(x)\mathcal{M}f_2(x).
\end{align*}
This completes the proof of \eqref{minorptmest} because $2^{\mu-10}\le C_0$.

\hfill

\section{Proof of Proposition \ref{initialpropo}}\label{Sec: Section 6}

Assume that $2^{\mu-10}>C_0$. Using \eqref{daubechewavelet}, we write
\begin{align*}
  \big(\Psi_{\mu-\gamma}\ast K^{*}_{\gamma}\big)^{\wedge}(\xxxi/2^{\gamma})&=\wh{\Psi}(2^{-\mu}\xxxi)\wh{K^{*}_{\gamma}}(2^{-\gamma}\xxxi)\\
  &=\sum_{\lambda\in\bbn_0}\sum_{\GGG\in \II^{\lambda}}\sum_{\kkk\in\bbz^2} b_{\GGG,\kkk}^{\lambda,\gamma,\mu}\vartheta_{G_1,k_1}^{\lambda}(\xi_1) \vartheta_{G_2,k_2}^{\lambda}(\xi_2)
\end{align*}
where
$$b_{\GGG,\kkk}^{\lambda,\gamma,\mu}:=\int_{\bbr^2} \wh{\Psi}(\xxxi/2^{\mu})\wh{K^{*}_{\gamma}}(2^{-\gamma}\xxxi)\Theta_{\GGG,\kkk}^{\lambda}(\xxxi)    \; d\xxxi.$$

We claim that for any $\frac{4}{3}<q\le 2$ and $0<\delta_0<\frac{1}{q'}$,
\begin{align}
  \big\Vert \big\{ b_{\GGG,\kkk}^{\lambda,\gamma,\mu}\big\}_{\kkk\in\bbz^2}\big\Vert_{\ell^{\infty}}&\lesssim 2^{-\delta_0\mu}2^{-\lambda (L+2)}\Vert \Omega\Vert_{L^q(\mathbb{S}^1)}\label{bdiscest1}\\
  \big\Vert \big\{ b_{\GGG,\kkk}^{\lambda,\gamma,\mu}\big\}_{\kkk\in\bbz^2}\big\Vert_{\ell^{q'}}&\lesssim 2^{-\lambda (\frac{2}{q}-1)}\Vert \Omega\Vert_{L^q(\mathbb{S}^1)} \label{bdiscest2}
\end{align}
uniformly in the parameters $\lambda,\gamma,\mu,\GGG$.
For this one, we first see that
$$\Big|  \partial^{\alpha}\Big( \wh{\Psi}(2^{-\mu}\vec{\cdot}\,)\wh{K^{*}_{\gamma}}(2^{-\gamma}\vec{\cdot}\,)\Big)(\xxxi)   \Big| \lesssim {2^{-\delta_0\mu}}\Vert \Omega\Vert_{L^q(\mathbb{S}^1)}\frac{1}{(1+2^{-\mu}|\xxxi|)^{M}}$$
and then \eqref{bdiscest1} follows immediately from \cite[Lemma 2.1]{Gr_He_Sl2020}.
To achieve \eqref{bdiscest2}, we apply \eqref{lqestimate} to obtain
$$ \big\Vert \big\{ b_{\GGG,\kkk}^{\lambda,\gamma,\mu}\big\}_{\kkk\in\bbz^2}\big\Vert_{\ell^{q'}}\lesssim 2^{-\lambda(1-\frac{2}{q'})}\big\Vert  \big(\Psi_{\mu-\gamma}\ast K^{*}_{\gamma}\big)^{\wedge}(2^{-\gamma}\vec{\cdot}\,)    \big\Vert_{L^{q'}(\bbr^2)}$$
and the right-hand side is controlled by
\begin{align*}
  2^{-\lambda(1-\frac{2}{q'})}2^{2\gamma}\big\Vert   \Psi_{\mu-\gamma}\ast K^{*}_{\gamma}(2^\gamma \, \cdot)  \big\Vert_{L^q(\bbr^2)}&=2^{-\lambda(\frac{2}{q}-1)}2^{2\gamma(1-\frac{1}{q})}\big\Vert   \Psi_{\mu-\gamma}\ast K^{*}_{\gamma} \big\Vert_{L^q(\bbr^2)}\\
  &\le 2^{-\lambda(\frac{2}{q}-1)}2^{2\gamma(1-\frac{1}{q})}\Vert \Psi_{\mu-\gamma}\Vert_{L^1(\bbr^2)} \Vert K^{*}_{\gamma}\Vert_{L^q(\bbr^2)}\\
  &\lesssim 2^{-\lambda(\frac{2}{q}-1)}2^{2\gamma(1-\frac{1}{q})}\big\Vert K \chi_{|\cdot|\sim 2^{\gamma}}\big\Vert_{L^q(\bbr^2)}\\
  &\lesssim 2^{-\lambda(\frac{2}{q}-1)}\Vert \Omega\Vert_{L^q(\mathbb{S}^1)}.
\end{align*}
This verifies \eqref{bdiscest2}.

Now, $\LL_{\Omega,\mu}^{\gamma}\big(f_1,f_2\big)(x)$ can be written as
\begin{equation}\label{lommugaex}
 \LL_{\Omega,\mu}^{\gamma}\big(f_1,f_2\big)(x)=\sum_{\lambda\in\bbn_0}\sum_{\GGG\in \II^{\lambda}}\sum_{\kkk\in\bbz^2: 2^{\lambda+\mu-2}\le |\kkk|\le 2^{\lambda+\mu+2}} b_{\GGG,\kkk}^{\lambda,\gamma,\mu}\Gamma_{G_1,k_1}^{\lambda,\gamma}f_1(x)\Gamma_{G_2,k_2}^{\lambda,\gamma}f_2(x)
 \end{equation}
where the operator $\Gamma_{G,k}^{\la,\ga}$ is defined via the Fourier transform by
\begin{equation*}
  \big( \Gamma_{G,k}^{\la,\ga}f\big)^{\wedge}(\xi):= \vartheta_{G,k}^{\la}(2^{-\gamma}\xi)\wh{f}(\xi).
\end{equation*} 
The sum over $\kkk$ was originally taken over all $\kkk\in\bbz^2$, but it could be replaced by the sum over $2^{\lambda+\mu-2}|\kkk|\le 2^{\lambda+\mu+2}$ because $b_{\GGG,\kkk}^{\lambda,\gamma,\mu}$ disappears unless $2^{\lambda+\mu-2}|\kkk|\le 2^{\lambda+\mu+2}$,
due to the compact support conditions of $\wh{\Psi}(2^{-\mu}\cdot)$ and $\Theta_{\GGG,\kkk}^{\lambda}$. Such materials were already discussed in \cite[Section 6]{Gr_He_Ho_Park2023}.
We additionally assume that $|k_1|\ge |k_2|$ as the opposite case follows by symmetry. Now, setting
$$\UU^{\lambda+\mu}:=\{\kkk\in\bbz^2: 2^{\lambda+\mu-2}\le |\kkk|\le 2^{\lambda+\mu+2}, \; |k_1|\le |k_2|\},$$
we claim that there exist $M>0$ and $\delta>0$ such that
\begin{equation}\label{propo1goal}
  \bigg\Vert \sup_{\gamma\in\bbz} \Big|  \sum_{\kkk\in\UU^{\lambda+\mu}}   b_{\GGG,\kkk}^{\lambda,\gamma,\mu}\Gamma_{G_1,k_1}^{\lambda,\gamma}f_1\Gamma_{G_2,k_2}^{\lambda,\gamma}f_2        \Big| \bigg\Vert_{L^1(\bbr)}\lesssim 2^{-\lambda M}{2^{-\delta \mu}}\Vert\Omega\Vert_{L^q(\mathbb{S}^1)}\Vert f_1\Vert_{L^2(\bbr)}\Vert f_2\Vert_{L^2(\bbr)},
\end{equation}
which clearly concludes \eqref{initialpropoest} in view of \eqref{lommugaex}.

To achieve \eqref{propo1goal}, we decompose $\UU^{\lambda+\mu}$ as 
\begin{equation}\label{udisdecom}
  \UU^{\lambda+\mu}=\UU_1^{\lambda+\mu} \;\dot{\cup}\; \UU_2^{\lambda+\mu}
\end{equation}
where
$$\UU_1^{\lambda+\mu}:=\big\{\kkk\in\UU^{\lambda+\mu}: |k_1|\ge 2C_0>|k_2| \big\},$$
$$\UU_2^{\lambda+\mu}:=\big\{\kkk\in\UU^{\lambda+\mu}: |k_1|\ge |k_2|\ge 2C_0 \big\},$$
as $|k_1|<2C_0$ implies $2^{\lambda+\mu-2}\le |\kkk|<2\sqrt{2}C_0$ and this is impossible.
Using \eqref{udisdecom},
we write
\begin{align*}
  &\bigg\Vert \sup_{\gamma\in\bbz}\Big|\sum_{\kkk\in\UU^{\lambda+\mu}} b_{\GGG,\kkk}^{\lambda,\gamma,\mu}\Gamma_{G_1,k_1}^{\lambda,\gamma}f_1\Gamma_{G_2,k_2}^{\lambda,\gamma}f_2        \Big|  \bigg\Vert_{L^1(\bbr)}\\
  &\le \bigg\Vert\sup_{\gamma\in\bbz} \Big|  \sum_{\kkk\in\UU_1^{\lambda+\mu}} b_{\GGG,\kkk}^{\lambda,\gamma,\mu}\Gamma_{G_1,k_1}^{\lambda,\gamma}f_1\Gamma_{G_2,k_2}^{\lambda,\gamma}f_2        \Big|\bigg\Vert_{L^1(\bbr)}\\
   &\quad + \bigg\Vert \sup_{\gamma\in\bbz} \Big|  \sum_{\kkk\in\UU_2^{\lambda+\mu}} b_{\GGG,\kkk}^{\lambda,\gamma,\mu}\Gamma_{G_1,k_1}^{\lambda,\gamma}f_1\Gamma_{G_2,k_2}^{\lambda,\gamma}f_2        \Big|\bigg\Vert_{L^1(\bbr)}\\
  &=:\mathfrak{I}_1+\mathfrak{I}_2
\end{align*}
and then will deal with the two terms in the right-hand side separately.
For this one, we observe several important estimates.

For each $G\in \{F,M\}$, $\lambda\in\bbn_0$, and $\gamma\in \bbz$,
the operator $\Gamma_{G,k}^{\la,\ga}$ satisfies

\begin{equation}\label{l2ovkest0}
  \big| \Gamma_{G,k}^{\lambda,\gamma}f(x)\big|  \lesssim 2^{\frac{\lambda}{2}}\mathcal{M}f(x)
\end{equation}
and
\begin{align}\label{l2ovkest}
  \bigg\Vert    \Big(   \sum_{k\in\bbz} \big| \Gamma_{G,k}^{\lambda,\gamma}f\big|^2    \Big)^{\frac{1}{2}} \bigg\Vert_{L^2(\bbr)}&=\bigg( \sum_{k\in\bbz}\big\Vert \Gamma_{G,k}^{\lambda,\gamma}f\big\Vert_{L^2(\bbr)}^2\bigg)^{\frac{1}{2}}\nonumber\\
  &=\bigg( \sum_{k\in\bbz}\big\Vert  \vartheta_{G,k}^{\lambda}(2^{-\gamma}\cdot)\wh{f}\,\big\Vert_{L^2(\bbr)}^2\bigg)^{\frac{1}{2}}\nonumber\\
  &=\bigg( \int_{\bbr} \big| \wh{f}(\xi)\big|^2\sum_{k\in\bbz}\big| \vartheta_{G,k}^{\lambda}(2^{-\gamma}\xi)\big|^2\bigg)^{\frac{1}{2}}\nonumber\\
  &\lesssim 2^{\frac{\lambda}{2}}\Vert f\Vert_{L^2(\bbr)}
\end{align}
where we applied Plancherel's identity and \eqref{lrcondition}.

For $\la \in\bbn_0$ satisfying $C_0\le 2^{\la+1}$, let
$$\WW^{\la}:=\big\{ k\in\bbz: 2C_0\le |k|\le 2^{\la+2}\big\}.$$
Then it is known in \cite[(39)]{Gr_He_Ho_Park2023} that
for $k\in \WW^{\la+\mu}$ with $C_0\le 2^{\la+\mu+1}$,
\begin{equation}\label{lgkest}
  L_{G,k}^{\la,\ga}f=L_{G,k}^{\la,\ga}f^{\la,\ga,{\mu}}   
\end{equation}
due to the support of $\vartheta_{G}$, 
where 
\begin{equation*}
  \wh{f^{\la,\ga,\mu}}(\xi):=\wh{f}(\xi)\chi_{c_02^{\gamma-\lambda}\le |\xi|\le 2^{\gamma+\mu+3}}
\end{equation*} for some $c_0>0$, depending on $C_0$.
In this case, Plancherel's identity deduces
\begin{equation}\label{pcidkey}
  \Big\Vert \Big(\sum_{\ga\in\bbz}\big|f^{\la,\ga,\mu} \big|^2 \Big)^{\frac{1}{2}}\Big\Vert_{L^2(\bbr)} \lesssim (\mu+\la)^{\frac{1}{2}}\Vert f\Vert_{L^2(\bbr)},
\end{equation}
See \cite[(40)]{Gr_He_Ho_Park2023} for more details.

Now, we return to estimating $\mathfrak{I}_1$ and $\mathfrak{I}_2$.
If $\kkk\in\UU_1^{\lambda+\mu}$, then \eqref{lgkest} implies
\begin{equation*}
L_{G,k}^{\la,\ga}f_1=L_{G,k}^{\la,\ga}f_1^{\la,\ga,{\mu}}.
\end{equation*}
Therefore, applying \eqref{disjointdecom1} and \eqref{l2ovkest0}
\begin{align*}
  &\Big|\sum_{\kkk\in\UU_1^{\lambda+\mu}}   b_{\GGG,\kkk}^{\lambda,\gamma,\mu}\Gamma_{G_1,k_1}^{\lambda,\gamma}f_1(x)\Gamma_{G_2,k_2}^{\lambda,\gamma}f_2(x)\Big|\\
 &\le  \sum_{k_2\in \mathcal{P}_2\UU_1^{\lambda+\mu}} \big| \Gamma_{G_2,k_2}^{\lambda,\gamma}f_2(x)\big|\Big| \sum_{k_1\in Col_{k_2}^{\UU_1^{\lambda+\mu},1}} b_{\GGG,\kkk}^{\lambda,\gamma,\mu}\Gamma_{G_1,k_1}^{\lambda,\gamma}f_1^{\lambda,\gamma,\mu}(x) \Big|\\
 &\lesssim 2^{\frac{\lambda}{2}}\mathcal{M}f_2(x) \sum_{k_2\in \mathcal{P}_2\UU_1^{\lambda+\mu}} \Big| \sum_{k_1\in Col_{k_2}^{\UU_1^{\lambda+\mu},1}} b_{\GGG,\kkk}^{\lambda,\gamma,\mu}\Gamma_{G_1,k_1}^{\lambda,\gamma}f_1^{\lambda,\gamma,\mu}(x) \Big|
\end{align*}
and this yields that
\begin{align*}
  \mathfrak{I}_1&\lesssim 2^{\frac{\lambda}{2}}\big\Vert \mathcal{M}f_2\big\Vert_{L^2(\bbr)}\bigg\Vert \sup_{\gamma\in\bbz}\sum_{k_2\in\mathcal{P}_2\UU_1^{\lambda+\mu}}\Big| \sum_{k_1\in Col_{k_2}^{\UU_1^{\lambda+\mu},1}} b_{\GGG,\kkk}^{\lambda,\gamma,\mu}\Gamma_{G_1,k_1}^{\lambda,\gamma}f_1^{\lambda,\gamma,\mu}  \Big|\bigg\Vert_{L^2(\bbr)}\\
  &\lesssim 2^{\frac{\lambda}{2}}\Vert f_2\Vert_{L^2(\bbr)}\bigg( \sum_{\gamma\in\bbz}\sum_{k_2\in \mathcal{P}_2\UU_1^{\lambda+\mu}}\Big\Vert \sum_{k_1\in Col_{k_2}^{\UU_1^{\lambda+\mu},1}}   b_{\GGG,\kkk}^{\lambda,\gamma,\mu}\Gamma_{G_1,k_1}^{\lambda,\gamma}f_1^{\lambda,\gamma,\mu}  \Big\Vert_{L^2(\bbr)}^2\bigg)^{\frac{1}{2}}
\end{align*}
where we applied Cauchy-Schwartz's inequality and the fact that 
$$\big| \mathcal{P}_2\UU_1^{\lambda+\mu}\big|\lesssim_{C_0} 1.$$
 Plancherel's identity and \eqref{lrcondition} deduces
\begin{align*}
  \Big\Vert \sum_{k_1\in Col_{k_2}^{\UU_1^{\lambda+\mu},1}} b_{\GGG,\kkk}^{\lambda,\gamma,\mu}\Gamma_{G_1,k_1}^{\lambda,\gamma}f_1^{\lambda,\gamma,\mu} \Big\Vert_{L^2(\bbr)}&\sim \Big\Vert \sum_{k_1\in Col_{k_2}^{\UU_1^{\lambda+\mu},1}} b_{\GGG,\kkk}^{\lambda,\gamma,\mu}\vartheta_{G_1,k_1}^{\lambda}(2^{-\gamma}\cdot)\wh{f_1^{\lambda,\gamma,\mu}}  \Big\Vert_{L^2(\bbr)}\\
  &\lesssim 2^{\frac{\lambda}{2}}\big\Vert \{b_{\GGG,\kkk}^{\lambda,\gamma,\mu}\}_{\kkk\in\bbz}\big\Vert_{\ell^{\infty}} \Big\Vert\wh{f_1^{\lambda,\gamma,\mu}}     \Big\Vert_{L^2(\bbr)}\\
   &=2^{\frac{\lambda}{2}}\big\Vert \{b_{\GGG,\kkk}^{\lambda,\gamma,\mu}\}_{\kkk\in\bbz}\big\Vert_{\ell^{\infty}} \Big\Vert f_1^{\lambda,\gamma,\mu}     \Big\Vert_{L^2(\bbr)}.
\end{align*}
This proves
\begin{align}
  \mathfrak{I}_1&\lesssim 2^{\lambda}\Vert f_2\Vert_{L^2(\bbr)} \Big( \sup_{\gamma\in\bbz}\big\Vert \{b_{\GGG,\kkk}^{\lambda,\gamma,\mu}\}_{\kkk\in\bbz}\big\Vert_{\ell^{\infty}} \Big)\bigg\Vert \Big(\sum_{\gamma\in\bbz} \big| f_1^{\lambda,\gamma,\mu}\big|^2 \Big)^{\frac{1}{2}}\bigg\Vert_{L^2(\bbr)}\nonumber\\
  &\lesssim 2^{\lambda}(\lambda+\mu)^{\frac{1}{2}} \sup_{\gamma\in\bbz}\big\Vert \{b_{\GGG,\kkk}^{\lambda,\gamma,\mu}\}_{\kkk\in\bbz}\big\Vert_{\ell^{\infty}}\Vert f_1\Vert_{L^2(\bbr)}  \Vert f_2\Vert_{L^2(\bbr)}\nonumber\\
  &\lesssim 2^{-\lambda(L+1)}2^{-\delta_0\mu}(\lambda+\mu)^{\frac{1}{2}}\Vert \Omega\Vert_{L^q(\mathbb{S}^1)}\Vert f_1\Vert_{L^2(\bbr)}\Vert f_2\Vert_{L^2(\bbr)} \label{j1mainest}
\end{align}
where the last two estimates follow from \eqref{pcidkey}, \eqref{bdiscest1}, and \eqref{bdiscest2}.

On the other hand, if  $\kkk\in\UU_2^{\lambda+\mu}$, then we have
\begin{equation*}
  L_{G,k}^{\la,\ga}f_1=L_{G,k}^{\la,\ga}f_1^{\la,\ga,{\mu}} ~\text{ and }~ L_{G,k}^{\la,\ga}f_2=L_{G,k}^{\la,\ga}f_2^{\la,\ga,{\mu}}
\end{equation*}
and thus we may replace $f_1$ and $f_2$ by $f_1^{\lambda,\gamma,\mu}$ and $f_2^{\lambda,\gamma,\mu}$, respectively in $\mathfrak{I}_2$.
For each $r\in\bbn$, let
$$\UU_{2,r}^{\lambda+\mu}:=\Big\{\kkk\in\UU_2^{\lambda+\mu}:2^{-r} \big\Vert \{b_{\GGG,\kkk}^{\lambda,\gamma,\mu}\}_{\kkk\in\bbz}\big\Vert_{\ell^{\infty}}<  \big| b_{\GGG,\kkk}^{\lambda,\gamma,\mu} \big|\le 2^{-r+1}\big\Vert\{b_{\GGG,\kkk}^{\lambda,\gamma,\mu}\}_{\kkk\in\bbz}\big\Vert_{\ell^{\infty}}\Big\}$$
so that
$$\UU_2^{\lambda+\mu}=\dot{\bigcup_{r\in\bbn}}\UU_{2,r}^{\lambda+\mu}.$$
We note that for each $r\in\bbn$,
$$2^{-r}\big\Vert \{b_{\GGG,\kkk}^{\lambda,\gamma,\mu}\}_{\kkk\in\bbz}\big\Vert_{\ell^{\infty}} \big| \UU_{2,r}^{\lambda+\mu} \big|^{\frac{1}{q'}}\le \Big( \sum_{\kkk\in \UU_{2,r}^{\lambda+\mu}} \big| b_{\GGG,\kkk}^{\lambda,\gamma,\mu}\big|^{q'}\Big)^{\frac{1}{q'}}\le \big\Vert \{b_{\GGG,\kkk}^{\lambda,\gamma,\mu}\}_{\kkk\in\bbz}\big\Vert_{\ell^{q'}},$$
which implies
\begin{equation}\label{uu2rlmcard}
  \big| \UU_{2,r}^{\lambda+\mu} \big|\le  2^{rq'}\frac{\big\Vert \{b_{\GGG,\kkk}^{\lambda,\gamma,\mu}\}_{\kkk\in\bbz}\big\Vert_{\ell^{q'}}^{q'}}{      \big\Vert \{b_{\GGG,\kkk}^{\lambda,\gamma,\mu}\}_{\kkk\in\bbz}\big\Vert_{\ell^{\infty}}^{q'}}.
\end{equation}
We further decompose $\UU_{2,r}^{\lambda+\mu}$ into
\begin{align*}
  \UU_{2,r}^{\lambda+\mu,1}&:= \Big\{\kkk\in\UU_{2,r}^{\lambda+\mu}: \big| Col_{k_2}^{\UU_{2,r}^{\lambda+\mu},1}\big|>\big| \UU_{2,r}^{\lambda+\mu}\big|^{\frac{1}{2}} \Big\}\\
  \UU_{2,r}^{\lambda+\mu,2}&:= \Big\{\kkk\in\UU_{2,r}^{\lambda+\mu}: \big| Col_{k_2}^{\UU_{2,r}^{\lambda+\mu},1}\big|\le \big| \UU_{2,r}^{\lambda+\mu}\big|^{\frac{1}{2}} \Big\}.
\end{align*}
Then we observe that
$$\big| \UU_{2,r}^{\lambda+\mu}\big|\ge \big| \UU_{2,r}^{\lambda+\mu,1}\big|\ge \big| \mathcal{P}_2\UU_{2,r}^{\lambda+\mu,1}\big|\big| \UU_{2,r}^{\lambda+\mu}\big|^{\frac{1}{2}},$$
which deduces
\begin{equation}\label{pro2u2rcard}
  \big| \mathcal{P}_2\UU_{2,r}^{\lambda+\mu,1}\big|\le \big| \UU_{2,r}^{\lambda+\mu}\big|^{\frac{1}{2}}.
\end{equation}
In addition, for each $k_2\in \mathcal{P}_2\UU_{2,r}^{\lambda+\mu,2}$
\begin{equation}\label{colk2u2rcard}
  \big| Col_{k_2}^{\UU_{2,r}^{\lambda+\mu,2},1}\big|\le \big| Col_{k_2}^{\UU_{2,r}^{\lambda+\mu},1}\big| \le \big| \UU_{2,r}^{\lambda+\mu}\big|^{\frac{1}{2}}.
\end{equation}

Now, we have
$$\mathfrak{I}_2\le \sum_{j=1}^{2}\sum_{r\in\bbn}\bigg\Vert  \sum_{\gamma\in\bbz}\Big|  \sum_{\kkk\in\UU_{2,r}^{\lambda+\mu,j}}   b_{\GGG,\kkk}^{\lambda,\gamma,\mu}\Gamma_{G_1,k_1}^{\lambda,\gamma}f_1\Gamma_{G_2,k_2}^{\lambda,\gamma}f_2\Big|\bigg\Vert_{L^1(\bbr)}$$
Then we claim that for each $j\in \{1,2\}$ and $r\in\bbn$,
\begin{align}\label{l1sgzbskiu}
  &\bigg\Vert  \sum_{\gamma\in\bbz}\Big|  \sum_{\kkk\in\UU_{2,r}^{\lambda+\mu,j}}   b_{\GGG,\kkk}^{\lambda,\gamma,\mu}\Gamma_{G_1,k_1}^{\lambda,\gamma}f_1\Gamma_{G_2,k_2}^{\lambda,\gamma}f_2      \Big|\bigg\Vert_{L^1(\bbr)}\nonumber\\
  &\lesssim 2^{\lambda}(\lambda+\mu)2^{-r(1-\frac{q'}{4})}\big\Vert \{b_{\GGG,\kkk}^{\lambda,\gamma,\mu}\}_{\kkk\in\bbz}\big\Vert_{\ell^{q'}}^{\frac{q'}{4}}\big\Vert \{b_{\GGG,\kkk}^{\lambda,\gamma,\mu}\}_{\kkk\in\bbz}\big\Vert_{\ell^{\infty}}^{1-\frac{q'}{4}}\Vert f_1\Vert_{L^2(\bbr)}\Vert f_2\Vert_{L^2(\bbr)}.
\end{align}
If the estimate \eqref{l1sgzbskiu} holds, then
we obtain, since $2\le q'<4$, 
\begin{align}
  \mathfrak{I}_2&\lesssim_q   2^{\lambda}(\lambda+\mu)\big\Vert \{b_{\GGG,\kkk}^{\lambda,\gamma,\mu}\}_{\kkk\in\bbz}\big\Vert_{\ell^{q'}}^{\frac{q'}{4}}\big\Vert \{b_{\GGG,\kkk}^{\lambda,\gamma,\mu}\}_{\kkk\in\bbz}\big\Vert_{\ell^{\infty}}^{1-\frac{q'}{4}} \Vert f_1\Vert_{L^2(\bbr)}\Vert f_2\Vert_{L^2(\bbr)}\nonumber\\
  &\lesssim 2^{-\lambda[(L+2)(1-\frac{q'}{4})+(\frac{2}{q}-1)\frac{q'}{4}-1]}2^{-\delta_0\mu(1-\frac{q'}{4})}(\lambda+\mu)\Vert \Omega\Vert_{L^q(\mathbb{S}^1)}\Vert f_1\Vert_{L^2(\bbr)}\Vert f_2\Vert_{L^2(\bbr)}\nonumber\\
  &=2^{-\lambda\frac{(3q-4)L+q-2}{4(q-1)}}2^{-\delta_0\mu\frac{3q-4}{4(q-1)}}(\lambda+\mu) \Vert \Omega\Vert_{L^q(\mathbb{S}^1)}\Vert f_1\Vert_{L^2(\bbr)}\Vert f_2\Vert_{L^2(\bbr)} \label{j2mainest}
\end{align}
where the estimates \eqref{bdiscest1} and \eqref{bdiscest2} are applied in the last inequality.
Accordingly, combining \eqref{j1mainest} and \eqref{j2mainest} together,
\eqref{propo1goal} follows from setting 
$$0<M<\frac{(3q-4)L+q-2}{4(q-1)}(<L+1) \quad \text{and} \quad 0<\delta<\delta_0\frac{3q-4}{4(q-1)}(<\delta_0).$$

Therefore, it remains to prove \eqref{l1sgzbskiu}.To begin with, consider the case $j=1$. By applying \eqref{lrcondition} and Cauchy-Schwartz's inequality,
\begin{align*}
  &\Big|  \sum_{\kkk\in\UU_{2,r}^{\lambda+\mu,1}}   b_{\GGG,\kkk}^{\lambda,\gamma,\mu}\Gamma_{G_1,k_1}^{\lambda,\gamma}f_1(x)\Gamma_{G_2,k_2}^{\lambda,\gamma}f_2(x)        \Big|\\
  &\le  \sum_{k_2\in \mathcal{P}_{2}\UU_{2,r}^{\lambda+\mu,1}} \big| \Gamma_{G_2,k_2}^{\lambda,\gamma}f_2^{\lambda,\gamma,\mu}(x)\big|\Big| \sum_{k_1\in Col_{k_2}^{\UU_{2,r}^{\lambda+\mu,1},1}} b_{\GGG,\kkk}^{\lambda,\gamma,\mu}\Gamma_{G_1,k_1}^{\lambda,\gamma}f_1^{\lambda,\gamma,\mu}(x) \Big|\\
  &\lesssim\bigg( \sum_{k_2\in \mathcal{P}_2\UU_{2,r}^{\lambda+\mu,1}}  \big| \Gamma_{G_2,k_2}^{\lambda,\gamma}f_2^{\lambda,\gamma,\mu}(x)\big|^2    \bigg)^{\frac{1}{2}}\\      &\qquad\times\bigg(\sum_{k_2\in \mathcal{P}_2\UU_{2,r}^{\lambda+\mu,1}} \Big| \sum_{k_1\in Col_{k_2}^{\UU_{2,r}^{\lambda+\mu,1},1}}      b_{\GGG,\kkk}^{\lambda,\gamma,\mu}\Gamma_{G_1,k_1}^{\lambda,\gamma}f_1^{\lambda,\gamma,\mu}(x) \Big|^2\bigg)^{\frac{1}{2}}.
\end{align*}
Therefore, by using Cauchy-Schwartz's inequality,
\begin{align}\label{u2rlm1est}
  &\bigg\Vert  \sum_{\gamma\in\bbz}\Big|  \sum_{\kkk\in\UU_{2,r}^{\lambda+\mu,1}}b_{\GGG,\kkk}^{\lambda,\gamma,\mu}\Gamma_{G_1,k_1}^{\lambda,\gamma}f_1\Gamma_{G_2,k_2}^{\lambda,\gamma}f_2 \Big|\bigg\Vert_{L^1(\bbr)}\nonumber\\
  &\le \bigg\Vert\bigg( \sum_{\gamma\in\bbz}\sum_{k_2\in \mathcal{P}_2\UU_{2,r}^{\lambda+\mu,1}}\big| \Gamma_{G_2,k_2}^{\lambda,\gamma}f_2^{\lambda,\gamma,\mu}\big|^2\bigg)^{\frac{1}{2}}\bigg\Vert_{L^2(\bbr)}\\
  &\qq\times \bigg\Vert \bigg(\sum_{\gamma\in\bbz}\sum_{k_2\in \mathcal{P}_2\UU_{2,r}^{\lambda+\mu,1}} \Big| \sum_{k_1\in Col_{k_2}^{\UU_{2,r}^{\lambda+\mu,1},1}}b_{\GGG,\kkk}^{\lambda,\gamma,\mu}\Gamma_{G_1,k_1}^{\lambda,\gamma}f_1^{\lambda,\gamma,\mu} \Big|^2\bigg)^{\frac{1}{2}}\bigg\Vert_{L^2(\bbr)}.\nonumber
\end{align}
The first $L^2$ norm (via \eqref{l2ovkest} and \eqref{pcidkey}) is dominated by a constant times
$$\lambda^{\frac{1}{2}}(\lambda+\mu)^{\frac{1}{2}}\Vert f_2\Vert_{L^2(\bbr)}.$$
The second $L^2$ norm is equal to
\begin{align}
  &\bigg( \sum_{\gamma\in\bbz}\sum_{k_2\in\mathcal{P}_2\UU_{2,r}^{\lambda+\mu,1}}\bigg\Vert      \sum_{k_1\in Col_{k_2}^{\UU_{2,r}^{\lambda+\mu,2},1}} b_{\GGG,\kkk}^{\lambda,\gamma,\mu}\Gamma_{G_1,k_1}^{\lambda,\gamma}f_1^{\lambda,\gamma,\mu}\bigg\Vert_{L^2(\bbr)}^2 \bigg)^{\frac{1}{2}}\nonumber\\
  &=\bigg( \sum_{\gamma\in\bbz}\sum_{k_2\in\mathcal{P}_2\UU_{2,r}^{\lambda+\mu,1}}\bigg\Vert    \sum_{k_1\in Col_{k_2}^{\UU_{2,r}^{\lambda+\mu,2},1}} b_{\GGG,\kkk}^{\lambda,\gamma,\mu}\vartheta_{G_1,k_1}^{\lambda}(2^{-\gamma}\cdot)\wh{f_1^{\lambda,\gamma,\mu}}      \bigg\Vert_{L^2(\bbr)}^2 \bigg)^{\frac{1}{2}}\nonumber\\
  &\le 2^{-r}\big\Vert \big\{b_{\GGG,\kkk}^{\lambda,\gamma,\mu} \big\}_{\kkk}\big\Vert_{\ell^{\infty}}\bigg( \int_{\bbr} \sum_{\gamma\in\bbz}\big|\wh{f_1^{\lambda,\gamma,\mu}}(\xi) \big|^2\sum_{k_2\in \mathcal{P}_2\UU_{2,r}^{\lambda+\mu,1}}\Big( \sum_{k_1\in Col_{k_2}^{\UU_{2,r}^{\lambda+\mu,1}}}\big| \vartheta_{G_1,k_1}^{\lambda}(2^{-\gamma}\xi)\big|^2 \Big)  \;   d\xi  \bigg)^{\frac{1}{2}}\label{2mrblinl2}
\end{align}
where we applied the disjoint support conditions of $\vartheta_{G_1,k_1}^{\lambda}(2^{-\gamma}\cdot)$. Now, we apply  \eqref{lrcondition}, \eqref{uu2rlmcard}, and \eqref{pro2u2rcard} to bound the last expression by
\begin{align*}
  &2^{\lambda}2^{-r}\big\Vert \big\{  b_{\GGG,\kkk}^{\lambda,\gamma,\mu} \big\}_{\kkk}\big\Vert_{\ell^{\infty}}  \bigg( \int_{\bbr} \sum_{\gamma\in\bbz}   \big| \wh{f_1^{\lambda,\gamma,\mu}}(\xi) \big|^2 \big| \mathcal{P}_2\UU_{2,r}^{\lambda+\mu,1}\big|   \;   d\xi  \bigg)^{\frac{1}{2}}\\
  &\lesssim 2^{\lambda}2^{-r(1-\frac{q'}{4})}\big\Vert \{b_{\GGG,\kkk}^{\lambda,\gamma,\mu}\}_{\kkk\in\bbz}\big\Vert_{\ell^{q'}}^{\frac{q'}{4}}\big\Vert \{b_{\GGG,\kkk}^{\lambda,\gamma,\mu}\}_{\kkk\in\bbz}\big\Vert_{\ell^{\infty}}^{1-\frac{q'}{4}},
\end{align*}
as desired.

On the other hand, when $j=2$, by symmetry, we could obtain
\begin{align*}
  &\bigg\Vert \sum_{\gamma\in\bbz}\Big|  \sum_{\kkk\in\UU_{2,r}^{\lambda+\mu,2}}   b_{\GGG,\kkk}^{\lambda,\gamma,\mu}\Gamma_{G_1,k_1}^{\lambda,\gamma}f_1\Gamma_{G_2,k_2}^{\lambda,\gamma}f_2        \Big|\bigg\Vert_{L^1(\bbr)}\\
  &\le \bigg\Vert  \bigg( \sum_{\gamma\in\bbz}\sum_{k_1\in \mathcal{P}_1\UU_{2,r}^{\lambda+\mu,2}}  \big| \Gamma_{G_1,k_1}^{\lambda,\gamma}f_1^{\lambda,\gamma,\mu}\big|^2\bigg)^{\frac{1}{2}}\bigg\Vert_{L^2(\bbr)}\\
  &\qq\times \bigg\Vert \bigg(  \sum_{\gamma\in\bbz}  \sum_{k_1\in \mathcal{P}_1\UU_{2,r}^{\lambda+\mu,2}} \Big| \sum_{k_2\in Col_{k_1}^{\UU_{2,r}^{\lambda+\mu,2},2}}b_{\GGG,\kkk}^{\lambda,\gamma,\mu}\Gamma_{G_2,k_2}^{\lambda,\gamma}f_2^{\lambda,\gamma,\mu} \Big|^2\bigg)^{\frac{1}{2}}    \bigg\Vert_{L^2(\bbr)},
\end{align*}
similar to \eqref{u2rlm1est}. The first term is clearly controlled by
$$\lambda^{\frac{1}{2}}(\lambda+\mu)^{\frac{1}{2}}\Vert f_1\Vert_{L^2(\bbr)}.$$
Moreover,  similar to \eqref{2mrblinl2}, we have
\begin{align*}
  &\bigg\Vert \bigg(\sum_{\gamma\in\bbz}\sum_{k_1\in \mathcal{P}_1\UU_{2,r}^{\lambda+\mu,2}} \Big| \sum_{k_2\in Col_{k_1}^{\UU_{2,r}^{\lambda+\mu,2},2}}b_{\GGG,\kkk}^{\lambda,\gamma,\mu}\Gamma_{G_2,k_2}^{\lambda,\gamma}f_2^{\lambda,\gamma,\mu} \Big|^2\bigg)^{\frac{1}{2}}    \bigg\Vert_{L^2(\bbr)}\\
  &\le 2^{-r}\big\Vert \big\{  b_{\GGG,\kkk}^{\lambda,\gamma,\mu} \big\}_{\kkk}\big\Vert_{\ell^{\infty}}  \bigg( \int_{\bbr} \sum_{\gamma\in\bbz}\big| \wh{f_2^{\lambda,\gamma,\mu}}(\xi) \big|^2\sum_{k_1\in \mathcal{P}_1\UU_{2,r}^{\lambda+\mu,2}}\Big( \sum_{k_2\in Col_{k_1}^{\UU_{2,r}^{\lambda+\mu,2}}}  \big| \vartheta_{G_2,k_2}^{\lambda}(2^{-\gamma}\xi)\big|^2 \Big)\;d\xi\bigg)^{\frac{1}{2}}.
\end{align*}
Now, we apply \eqref{disjointdecom1} again to obtain
\begin{align*}
  \sum_{k_1\in \mathcal{P}_1\UU_{2,r}^{\lambda+\mu,2}}\Big( \sum_{k_2\in Col_{k_1}^{\UU_{2,r}^{\lambda+\mu,2}}}\big| \vartheta_{G_2,k_2}^{\lambda}(2^{-\gamma}\xi)\big|^2 \Big) &=\sum_{k_2\in \mathcal{P}_2\UU_{2,r}^{\lambda+\mu,2}}\big| \vartheta_{G_2,k_2}^{\lambda}(2^{-\gamma}\xi)\big|^2\Big(\sum_{k_1\in Col_{k_2}^{\UU_{2,r}^{\lambda+\mu},1}}1 \Big)\\
  &\le 2^{\lambda}2^{\frac{rq'}{2}}\frac{\big\Vert \{b_{\GGG,\kkk}^{\lambda,\gamma,\mu}\}_{\kkk\in\bbz}\big\Vert_{\ell^{q'}}^{\frac{q'}{2}}}{\big\Vert\{b_{\GGG,\kkk}^{\lambda,\gamma,\mu}\}_{\kkk\in\bbz}\big\Vert_{\ell^{\infty}}^{\frac{q'}{2}}}
\end{align*}
where the last inequality follows from \eqref{colk2u2rcard}, \eqref{uu2rlmcard}, and \eqref{lrcondition}.
This proves
\begin{align*}
  &\bigg\Vert \bigg(  \sum_{\gamma\in\bbz}  \sum_{k_1\in \mathcal{P}_1\UU_{2,r}^{\lambda+\mu,2}} \Big| \sum_{k_2\in Col_{k_1}^{\UU_{2,r}^{\lambda+\mu,2},2}}b_{\GGG,\kkk}^{\lambda,\gamma,\mu}\Gamma_{G_2,k_2}^{\lambda,\gamma}f_2^{\lambda,\gamma,\mu} \Big|^2\bigg)^{\frac{1}{2}}  \bigg\Vert_{L^2(\bbr)}\\
  &\lesssim 2^{\frac{\lambda}{2}}2^{-r(1-\frac{q'}{4})}\big\Vert \{b_{\GGG,\kkk}^{\lambda,\gamma,\mu}\}_{\kkk\in\bbz}\big\Vert_{\ell^{q'}}^{\frac{q'}{4}}\big\Vert \{b_{\GGG,\kkk}^{\lambda,\gamma,\mu}\}_{\kkk\in\bbz}\big\Vert_{\ell^{\infty}}^{1-\frac{q'}{4}}.
\end{align*}
This completes the proof of \eqref{l1sgzbskiu}.

\hfill

\section{Proof of Proposition \ref{PropoHLS}}\label{Sec: Section 7}

Without loss of generality, we may assume 
$$\Vert f_1\Vert_{L^1(\bbr)}=\Vert f_2\Vert_{L^1(\bbr)}=\Vert \Omega\Vert_{L^1(\mathbb{S}^1,u_1)}=1$$
and then we will prove that for any $\epsilon>0$ and for any $\lambda>0$
\begin{equation}\label{explmainest}
  \bigg|\Big\{ x\in \bbr: \sup_{\gamma\in\bbz}\Big| \LL_{\Omega,\mu}^{\gamma}(f_1,f_2)(x)\Big|>\lambda\Big\} \bigg|\lesssim_{\epsilon} 2^{\epsilon \mu}\lambda^{-\frac{1}{2}}.
\end{equation}
We apply a Calder\'on-Zygmund decomposition of $f_1$ and $f_2$ at height $\lambda^{\frac{1}{2}}$, as in the proof of \cite[Proposition 4]{Ho_La_Sl_submitted}. To be specific,
we write
\begin{equation}\label{CZdecompff12}
  f_1=g_1+\sum_{Q\in\mathcal{A}}b_{1,Q}\qq \text{ and }\qq f_2=g_2+\sum_{R\in \mathcal{B}}b_{2,R}
\end{equation}
where
\begin{enumerate}
\item $\mathcal{A}$ and $\mathcal{B}$ are collections of non-overlapping dyadic intervals such that
$$\Big| \bigcup_{Q\in\mathcal{A}}Q\Big|\le \lambda^{-\frac{1}{2}}\qq \text{ and }\qq \Big| \bigcup_{R\in\mathcal{B}}R\Big|\le \lambda^{-\frac{1}{2}}$$
\item $b_{1,Q}$ and $b_{2,R}$ are supported in $Q$ and $R$, respectively
\item $$\int_{\bbr} b_{1,Q}(y)\; dy=\int_{\bbr} b_{2,R}(y)\; dy=0$$
\item $$\Vert b_{1,Q}\Vert_{L^1(\bbr)}\le 4\lambda^{\frac{1}{2}}|Q| \qq \text{ and }\qq \Vert b_{2,R}\Vert_{L^1(\bbr)}\le 4 \lambda^{\frac{1}{2}}|R|$$
\item For any $1\le r\le \infty$
$$\Vert g_1\Vert_{L^r(\bbr)},  \Vert g_2\Vert_{L^r(\bbr)} \le 2 \lambda^{\frac{1}{2}(1-\frac{1}{r})}.$$
\end{enumerate}
Moreover, for given $Q\in \mathcal{A}$ we also define
\begin{equation}\label{CZdecompff12:2}
  f_{2,Q}:=g_2+\sum_{R\in\mathcal{B}:\ell(R)\le \ell(Q)}b_{2,R}
\end{equation}
and similarly for $R\in \mathcal{B}$ we set
\begin{equation}\label{CZdecompff12:3}
  f_{1,R}:=g_1+\sum_{Q\in\mathcal{A}:\ell(Q)< \ell(R)}b_{1,Q}.
\end{equation}

Then the left-hand side of \eqref{explmainest} is bounded by the sum of
\begin{align*}
  \mathcal{I}_1&:=\bigg|\Big\{ x\in \bbr: \sup_{\gamma\in\bbz}\Big| \LL_{\Omega,\mu}^{\gamma}\big(g_1,g_2\big)(x)\Big|>\frac{\lambda}{3}\Big\} \bigg|\\
  \mathcal{I}_2&:=\bigg|\Big\{ x\in \bbr: \sum_{Q\in\mathcal{A}}\sup_{\gamma\in\bbz}\Big| \LL_{\Omega,\mu}^{\gamma}\big(b_{1,Q},f_{2,Q}\big)(x)\Big|>\frac{\lambda}{3}\Big\} \bigg|\\
  \mathcal{I}_3&:=\bigg|\Big\{ x\in \bbr: \sum_{R\in\mathcal{B}}\sup_{\gamma\in\bbz}\Big| \LL_{\Omega,\mu}^{\gamma}\big(f_{1,R},b_{2,R}\big)(x)\Big|>\frac{\lambda}{3}\Big\} \bigg|
\end{align*}

By applying Chebyshev's inequality and Proposition \ref{muge0tcpg1} with $L^1(\mathbb \S^1, u_1) \hookrightarrow L^{1}(\mathbb \S^1)$,
\begin{align*}
  \mathcal{I}_1\lesssim \lambda^{-2}\Big\Vert\sup_{\gamma\in\bbz}\big| \LL_{\Omega,\mu}^{\gamma}(g_1,g_2)\big| \Big\Vert_{L^2(\bbr)}^2\lesssim \lambda^{-2}\Vert \Omega\Vert_{L^1(\mathbb{S}^1)}^2\Vert g_1\Vert_{L^4(\bbr)}^2\Vert g_2\Vert_{L^4(\bbr)}^2\lesssim \lambda^{-\frac{1}{2}},
\end{align*}
as desired.

On the other hand, 
$$\mathcal{I}_2\le \Big| \bigcup_{Q\in\mathcal{A}}Q^*\Big|+ \bigg|\Big\{ x\in \Big(\bigcup_{Q\in\mathcal{A}}Q^*\Big)^c: \sum_{Q\in\mathcal{A}}\sup_{\gamma\in\bbz}\Big| \LL_{\Omega,\mu}^{\gamma}\big(b_{1,Q},f_{2,Q}\big)(x)\Big|>\frac{\lambda}{3}\Big\} \bigg|$$
where $Q^*$ is the concentric dilation of $Q$ with $\ell(Q^*)=10\ell(Q)$, and clearly,
$$\Big| \bigcup_{Q\in\mathcal{A}}Q^*\Big|\lesssim \lambda^{-\frac{1}{2}}.$$
Moreover, by applying Chebyshev's inequality and the embedding $\ell^1\hookrightarrow \ell^{\infty}$, we have
\begin{align*}
  &\bigg|\Big\{ x\in \Big(\bigcup_{Q\in\mathcal{A}}Q^*\Big)^c: \sum_{Q\in\mathcal{A}}\sup_{\gamma\in\bbz}\Big| \LL_{\Omega,\mu}^{\gamma}\big(b_{1,Q},f_{2,Q}\big)(x)\Big|>\frac{\lambda}{3}\Big\} \bigg|\\
  &\lesssim \lambda^{-1}\sum_{Q\in\mathcal{A}}\sum_{\gamma\in\bbz}\int_{(Q^*)^c}\Big| \LL_{\Omega,\mu}^{\gamma}\big(b_{1,Q},f_{2,Q}\big)(x)\Big| \; dx\\
  &=\lambda^{-1}\sum_{Q\in \mathcal{A}}\sum_{\gamma\in\bbz:2^{-\gamma}\ell(Q)> 1}\cdots \; + \lambda^{-1}\sum_{Q\in \mathcal{A}}\sum_{\gamma\in\bbz:2^{-\gamma}\ell(Q)\le  1}\cdots\\
  &=:\mathfrak{K}_1+\mathfrak{K}_2.
\end{align*}
Now we recall, following \cite[page 2267]{Do_Sl2024} and \cite[(3.6)]{Ho_La_Sl_submitted}, that
$$\big| \Psi(\zzz)\big|\lesssim\sum_{k=1}^{\infty}2^{-5k}\chi_{|\zzz|\le 2^k}(\zzz)$$
and thus
\begin{align}\label{pastkgtyest}
  \big|\Psi_{\mu-\gamma}\ast K^{*}_{\gamma}(\yyy)\big|& \lesssim \int_{2^{\gamma}\le |\zzz|\le 2^{\gamma+1}}\big| \Psi_{\mu-\gamma}(\yyy-\zzz)\big| |K(\zzz)|\; d\zzz \nonumber\\
  &\lesssim 2^{2\mu -2\gamma}\sum_{k=1}^{\infty}2^{-5k}\int_{1\le |\zzz|\le 2} \big| \Omega(\zzz')\big|\chi_{|\yyy-2^{\gamma}\zzz|\le 2^{k-\mu+\gamma}}\; d\zzz,
\end{align}
where \eqref{kbdptest} and a change of variables are applied in the estimates.
This yields that
\begin{align*}
  \mathfrak{K}_1&\lesssim \lambda^{-1}\sum_{Q\in\mathcal{A}}\sum_{\gamma\in\bbz: 2^{-\gamma}\ell(Q)>1}\int_{(Q^*)^c}\big| \LL_{\Omega,\mu}^{\gamma}\big(b_{1,Q},f_{2,Q}\big)(x)\big| \; dx\\
  &\lesssim \lambda^{-1}\sum_{Q\in\mathcal{A}}\sum_{\gamma\in\bbz: 2^{-\gamma}\ell(Q)>1}2^{2\mu-2\gamma}\sum_{k=1}^{\infty}2^{-5k}\int_{\bbr^2}\big| b_{1,Q}(y_1)f_{2,Q}(y_2)\big|\\
  &\qq\qq\times \bigg( \int_{1\le |\zzz|\le 2}\big|\Omega(\zzz') \big|\int_{(Q^*)^c}\chi_{|(x-y_1-2^{\gamma}z_1,x-z_2-2^{\gamma}z_2)|\le 2^{k-\mu+\gamma}}(x)\, dx\;d\zzz\bigg)\;d\yyy.
\end{align*}
We note that for $x\in (Q^*)^c$, $y_1\in Q$, $1\le |\zzz|\le 2$, $2^{-\gamma}\ell(Q)>1$, and $|(x-y_1-2^{\gamma}z_1,x-y_2-z_2)|\le 2^{k-\mu+\gamma}$,
$$\ell(Q)\le |x-y_1-2^{\gamma}z_1|\le 2^{k-\mu+\gamma}\le 2^{k+\gamma}$$
so that the sum over $k\in\bbn$ can be replaced by the sum over $k: 2^k\ge 2^{-\gamma}\ell(Q)$. 
In addition, 
\begin{align*}
 |y_1-y_2|\le |x-y_1-2^{\gamma}z_1|+ |x-y_2-2^{\gamma}z_2|+2^{\gamma}|z_1-z_2|\le 2^{k-\mu+\gamma+1}+2^{\gamma+1}\le 2^{k+2}\ell(Q),
\end{align*}
which yields that
$$y_2\in I(c_Q,2^{k+4}\ell(Q))$$
where $I(c_Q,2^{k+4}\ell(Q))$ denotes the concentric dilate of $Q$ whose length is $2^{k+4}\ell(Q)$.
We also recall from \cite{Ho_La_Sl_submitted}  that
\begin{align}\label{intomintosc}
  \int_{1\le |\zzz|\le 2}\big|\Omega(\zzz') \big|\int_{(Q^*)^c}\chi_{|(x-y_1-2^{\gamma}z_1,x-z_2-2^{\gamma}z_2)|\le 2^{k-\mu+\gamma}}(x)\, dx \;d\zzz\lesssim 2^{\gamma}2^{-2\mu +2k}
\end{align}
and
$$\Vert f_{2,Q}\Vert_{L^1(I(c_Q,2^{k+4}\ell(Q)))}\lesssim \lambda^{\frac{1}{2}}2^k\ell(Q).$$
See Lemma 5 and the inequality displayed after (3.8) in \cite{Ho_La_Sl_submitted} for more details.
Combining all together, we have
\begin{align*}
  \mathfrak{K}_1&\lesssim \lambda^{-1}\sum_{Q\in\mathcal{A}}\big\Vert b_{1,Q}\big\Vert_{L^1(Q)}\sum_{\gamma\in\bbz:2^{-\gamma}\ell(Q)>1}2^{-\gamma}\sum_{k\in\bbn:2^k\ge 2^{-\gamma}\ell(Q)}2^{-3k}\big\Vert f_{2,Q}\big\Vert_{L^1(I(c_Q,2^{k+4}\ell(Q)))}\\
  &\lesssim \lambda^{-\frac{1}{2}}\sum_{Q\in\mathcal{A}}\big\Vert b_{1,Q}\big\Vert_{L^1(Q)}\sum_{\gamma\in\bbz:2^{-\gamma}\ell(Q)>1}2^{-\gamma}\ell(Q)\sum_{k\in\bbn: 2^k \ge 2^{-\gamma}\ell(Q)}2^{-2k}\\
  &\lesssim \lambda^{-\frac{1}{2}}.
\end{align*}

To estimate $\mathfrak{K}_2$, we recall from \cite{Ho_La_Sl_submitted} that for any $0<\epsilon <1$
$$\big| \Psi(z_1,z_2)-\Psi(z,z_2)\big|\lesssim_{\epsilon} |z_1-z|^{\epsilon}\sum_{k=1}^{\infty}2^{-4k}\Big( \chi_{|(z_1,z_2)|\le 2^k}+\chi_{|(z,z_2)|\le 2^k}\Big)$$
and thus for each $\wt{y}\in\bbr$,
\begin{align*}
  \big|\Psi_{\mu-\gamma}\ast K^{*}_{\gamma}(\yyy)\big|& \lesssim \int_{2^{\gamma}\le |\zzz|\le 2^{\gamma+1}}\big| \Psi_{\mu-\gamma}(y_1-z_1,y_2-z_2)-\Psi_{\mu-\gamma}(\wt{y}-z_1,y_2-z_2)\big| |K(\zzz)|\; d\zzz\\
  &\lesssim 2^{2\mu -2\gamma}\big(2^{\mu-\gamma}|y_1-\wt{y}| \big)^{\epsilon}\sum_{k=1}^{\infty}2^{-4k}\int_{1\le |\zzz|\le 2} \big| \Omega(\zzz')\big| \\
  &\qq\qq\times\Big( \chi_{|(y_1-2^{\gamma}z_1,y_2-2^{\gamma}z_2)| \le 2^{k-\mu+\gamma}}+\chi_{|(\wt{y}-2^{\gamma}z_1,y_2-2^{\gamma}z_2)| \le 2^{k-\mu+\gamma}}  \Big)\; d\zzz,
\end{align*}
analogous to \eqref{pastkgtyest}.
Accordingly, we have
\begin{align*}
  \mathfrak{K}_2&\lesssim \lambda^{-1}\sum_{Q\in\mathcal{A}}\sum_{\gamma\in\bbz: 2^{-\gamma}\ell(Q)\le 1}2^{2\mu-2\gamma}\big(2^{\mu-\gamma}\ell(Q)\big)^{\epsilon}\sum_{k=1}^{\infty}2^{-4k}\int_{\bbr^2} \big| b_{1,Q}(y_1)f_{2,Q}(y_2)\big|\\
  &\qq\qq\times \bigg( \int_{1\le |\zzz|\le 2}  |\Omega(\zzz')|\int_{(Q^*)^c}  \chi_{|(x-y_1-2^{\gamma}z_1,x-y_2-2^{\gamma}z_2)|\le 2^{k-\mu+\gamma}}(x)    \; dx    \; d\zzz\bigg)\; d\yyy\\
  &\qq+\lambda^{-1}\sum_{Q\in\mathcal{A}}\sum_{\gamma\in\bbz: 2^{-\gamma}\ell(Q)\le 1}2^{2\mu-2\gamma}\big(2^{\mu-\gamma}\ell(Q)\big)^{\epsilon}\sum_{k=1}^{\infty}2^{-4k}\int_{\bbr^2} \big| b_{1,Q}(y_1)f_{2,Q}(y_2)\big|\\
  &\qq\qq\qq\times \bigg( \int_{1\le |\zzz|\le 2}  |\Omega(\zzz')|\int_{(Q^*)^c}  \chi_{|(x-c_Q-2^{\gamma}z_1,x-y_2-2^{\gamma}z_2)|\le 2^{k-\mu+\gamma}}(x)    \; dx    \; d\zzz\bigg)\; d\yyy\\
  &=:\mathfrak{K}_2^1+\mathfrak{K}_2^2.
\end{align*}
We will only consider $\mathfrak{K}_2^1$ as $\mathfrak{K}_2^2$ can be handled similarly.

For $x\in (Q^*)^c$, $y_1\in Q$, $1\le |\zzz|\le 2$, $2^{-\gamma}\ell(Q)\le 1$, and $|(x-y_1-2^{\gamma}z_1,x-y_2-2^{\gamma}z_2)|\le 2^{k-\mu+\gamma}$,
we have
$$|y_1-y_2|\le |x-y_1+2^{\gamma}z_1|+|x-y_2-2^{\gamma}z_2|+2^{\gamma}|z_1-z_2|\le 2^{k-\mu+\gamma+1}+2^{\gamma+1}\le 2^{k+\gamma+2}$$
and thus
$$y_2\in I(c_Q,2^{k+\gamma+4})$$
where we recall that $I(c_Q,2^{k+\gamma+4})$ stands for the interval centered at $c_Q$ whose length is $2^{k+\gamma+4}$.
Therefore, using \eqref{intomintosc},
\begin{align*}
  \mathfrak{K}_2^1&\lesssim \lambda^{-1}\sum_{Q\in\mathcal{A}}\big\Vert b_{1,Q} \big\Vert_{L^1(Q)} \sum_{\gamma\in\bbz: 2^{-\gamma}\ell(Q)\le 1}2^{-\gamma}\big(2^{\mu-\gamma}\ell(Q)\big)^{\epsilon}\sum_{k=1}^{\infty}2^{-2k} \big\Vert f_{2,Q}\big\Vert_{L^1(I(c_Q,2^{k+\gamma+4}))}.
\end{align*}
It is known in \cite[(3.11)]{Ho_La_Sl_submitted} that
$$\big\Vert f_{2,Q}\big\Vert_{L^1(I(c_Q,2^{k+\gamma+4}))}\lesssim \lambda^{\frac{1}{2}}2^{\gamma+k}$$
and this yields that
\begin{align*}
  \mathfrak{K}_2^1&\lesssim \lambda^{-\frac{1}{2}}\sum_{Q\in\mathcal{A}}\big\Vert b_{1,Q} \big\Vert_{L^1(Q)} \sum_{\gamma\in\bbz: 2^{-\gamma}\ell(Q)\le 1}\big(2^{\mu-\gamma}\ell(Q)\big)^{\epsilon}\\
  &\lesssim_{\epsilon} 2^{\epsilon \mu}\sum_{Q\in\mathcal{A}}|Q|\lesssim 2^{\epsilon \mu} \lambda^{-\frac{1}{2}},
\end{align*}
which deduces
$$\II_2\lesssim_{\epsilon} 2^{\epsilon \mu}\lambda^{-\frac{1}{2}}.$$

In a symmetric way, we can also establish
$$\II_3\lesssim_{\epsilon} 2^{\epsilon \mu}\lambda^{-\frac{1}{2}},$$
which completes the proof of \eqref{explmainest}.

\hfill

\section{Proof of Proposition \ref{prop:intmaxE}} \label{Sec: Section 8}

Suppose that $1<q<\infty$ and
$$\Vert f_1\Vert_{L^{p_1}(\bbr)}=\Vert f_2\Vert_{L^{p_2}(\bbr)}=\Vert \Omega\Vert_{L^{q}(\mathbb{S}^{1})}=1.$$
Then it is sufficient to show the existence of ${\delta}>0$ for which
\begin{equation}\label{mainkeyest33}
  \Big\|\sup_{\gamma\in\bbz}\big| \LL_{\Omega,\mu}^{\gamma}\big(f_1,f_2\big)\big| \Big\|_{L^{p}(\bbr) } \lesssim_{\delta}  2^{-\delta \mu}, \q \mu\in\bbz ~\text{with}~2^{\mu-10}>C_0 .
\end{equation}
For this one, we first decompose the circle $\mathbb{S}^{1}$ as
$$ \mathbb{S}^{1}=\dot{\bigcup_{l\in\bbn_0}}D^l$$
where
$$D^{l}:=\begin{cases} \big\{\theta\in\mathbb{S}^{1}: |\Omega(\theta)|\le 1\big\} & \text{ if } l=0\\
\big\{\theta\in\mathbb{S}^{1}:  2^{l-1}< |\Omega(\theta)|\le 2^l    \big\}  & \text{ if } l\ge 1 \end{cases}, $$
and write
\begin{equation*}
  \Omega(\theta)=\Omega(\theta)-\int_{\bbs^{1}}\Omega(\zeta)     \,d\sigma(\zeta) = \sum_{l=0}^{\infty}\bigg( \Omega(\theta)\chi_{D^{l}}(\theta)-\int_{D^{l}} \Omega(\zeta) \,d\sigma(\zeta)\bigg)=:\sum_{l=0}^{\infty}\Omega^{l}(\theta).
\end{equation*}
Then the left-hand side of \eqref{mainkeyest33} is bounded by
\begin{equation}\label{lastdecomex}
  \sum_{l\in\bbn_0} \Big\|   \sup_{\gamma\in\bbz}\big| \LL_{\Omega^l,\mu}^{\gamma}\big(f_1,f_2\big)\big|  \Big\|_{L^{p}(\bbr)}.
\end{equation}
We note that each $\Omega^{l}$ satisfies the vanishing moment condition 
$$\int_{\bbs^{1}}\Omega^{l}(\theta)\;d\sigma(\theta)=0$$
and thus we can apply  \eqref{generalrgestimp} and Proposition \ref{muge0tcpg1}  to $\Omega^l$ instead of $\Omega$.
Obviously,
\begin{equation*}
  \Vert\Omega^{l}\Vert_{L^{\infty}(\bbs^{1})}\le 2^{l+1}
\end{equation*}	
and thus \eqref{generalrgestimp} yields that for some $\delta_0>0$
\begin{equation}\label{commonest}
  \Big\Vert\sup_{\gamma\in\bbz}\big| \LL_{\Omega^l,\mu}^{\gamma}(f_1,f_2)\big|  \Big\Vert_{L^p(\bbr)}\lesssim_{\delta_0,C_0} 2^{-\delta_0\mu}\big\Vert \Omega^l\big\Vert_{L^{\infty}(\mathbb{S}^{1})}\lesssim 2^{-\delta_0 \mu}2^l
\end{equation}
Moreover, we see
\begin{equation*}
  \big\|\Omega^l\big\|_{L^1(\bbs^{1})} \le2\int_{D^l}\big|\Omega(\theta)\big| \;d\sigma(\theta)\lesssim_q 2^{-l(q-1)}\int_{D^l}\big| \Omega(\theta)\big|^{q}\;d\sigma(\theta)\le 2^{-lq}
\end{equation*}
and thus Proposition \ref{muge0tcpg1} deduces
\begin{equation}\label{commonestpge1}
  \Big\|\sup_{\gamma\in\bbz}\big| \LL_{\Om^l,\mu}^{\gamma}(f_1,f_2 )\big| \Big\|_{L^{p}(\bbr)}\lesssim  \big\|\Om^l\big\|_{L^1(\mathbb S^{1})}\lesssim 2^{-lq}.
\end{equation} 
We choose $\frac{1}{1+q}<\eta<1$, or consequently,
$$\eta q-(1-\eta)>0,$$
and by averaging \eqref{commonest} and \eqref{commonestpge1}, we obtain
\begin{align*}
  \Big\Vert\sup_{\gamma\in\bbz}\big|\LL_{\Omega^l,\mu}^{\gamma}(f_1,f_2)\big| \Big\Vert_{L^p(\bbr)}&\lesssim  \big(  2^{-lq} \big)^{\eta}\big( 2^{-\delta_0 \mu}2^l\big)^{1-\eta}=2^{-\delta_0(1-\eta)\mu} 2^{-l(\eta q-(1-\eta))}.
\end{align*}
Clearly, the right-hand side is summable over $l\in\bbn_0$ and thus  \eqref{lastdecomex} is dominated by a constant times
$$2^{-\delta_0(1-\eta)\mu}\Big(\sum_{l\in\bbn_0} 2^{-l(\eta q-(1-\eta) )}\Big)\lesssim_{\delta,C_0} 2^{-\delta \mu}$$
for some $0<\delta<\delta_0 (1-\eta)$, when $2^{\mu-10}>C_0$,
as desired.

\hfill

\section{Proof of Proposition \ref{prop:intmaxE1}}\label{Sec: Section 9}

Let $1<p_1,p_2<\infty$ and $\frac{1}{2}<p\le  1$ with $\frac{1}{p}=\frac{1}{p_1}+\frac{1}{p_2}$, and $q>1$.  In this case, we apply an interpolation method that is used in the proof of \eqref{lsharpommuest}, where we notice that   only the case $\frac{1}{2}<p<1$ was treated in the proof of \eqref{lsharpommuest}, but the same arguments are also valid for $p=1$.
In particular, as described in Figure \ref{fig:figure2}, we choose $\eta>0$ sufficiently small so that the point $A_1=(\frac{1}{p_1},\frac{1}{p_2})$ is located inside a non-trivial triangle whose vertices are $A_2=(1,1), A_3=(1-2\eta,\eta)$, and $A_4=(\eta,1-2\eta)$. According to Propositions \ref{PropoHLS} and \ref{prop:intmaxE}, for any $\epsilon>0$ and for some $\wt{\delta}>0$,
\begin{align*}
  \Big\Vert \sup_{\gamma\in\bbz}\big|\LL_{\Omega,\mu}^{\gamma}(f_1,f_2) \big|\Big\Vert_{ L^{\frac{1}{2},\infty}(\R)}&\lesssim_{\epsilon} 2^{\epsilon \mu}\Vert\Omega\Vert_{L^1(\S^{1}, u_1)}\Vert f_1\Vert_{L^{1}(\bbr)}\Vert f_2\Vert_{L^{1}(\bbr)},\\
  \Big\Vert \sup_{\gamma\in\bbz}\big|\LL_{\Omega,\mu}^{\gamma}(f_1,f_2) \big|\Big\Vert_{ L^{\frac{1}{1-\eta},\infty}(\R)}&\lesssim_{\wt{\delta}} 2^{-\wt{\delta} \mu}\Vert\Omega\Vert_{L^q(\S^{1})}\Vert f_1\Vert_{L^{\frac{1}{1-2\eta}}(\bbr)}\Vert f_2\Vert_{L^{\frac{1}{\eta}}(\bbr)},\\
  \Big\Vert \sup_{\gamma\in\bbz}\big|\LL_{\Omega,\mu}^{\gamma}(f_1,f_2) \big|\Big\Vert_{ L^{\frac{1}{1-\eta},\infty}(\R)}&\lesssim_{\wt{\delta}} 2^{-\wt{\delta} \mu}\Vert\Omega\Vert_{L^q(\S^{1})}\Vert f_1\Vert_{L^{\frac{1}{\eta}}(\bbr)}\Vert f_2\Vert_{L^{\frac{1}{1-2\eta}}(\bbr)}.
\end{align*}
Now, using embeddings $L^q(\S^{1},u_A^q)\hookrightarrow L^{q}(\S^{1})$ and $L^q(\S^{1},u_A^q)\hookrightarrow L^1(\S^{1},u_1)$ for $A>\frac{1}{q}$ and Lemma \ref{lem:standint} with \eqref{eq:main eq 3}, \eqref{eq:main eq 4}, and \eqref{eq:main eq 5},
it follows that
\begin{equation*}
  \Big\Vert \sup_{\gamma\in\bbz}\big|\LL_{\Omega,\mu}^{\gamma}(f_1,f_2) \big| \Big\Vert_{L^p(\bbr)} \lesssim 2^{-\mu(\wt{\delta} \theta_1+\wt{\delta} \theta_2-\epsilon \theta_0)     }\Vert\Omega\Vert_{L^q(\S^{1},u_A^q)}\Vert f_1\Vert_{L^{p_1}(\bbr)}\Vert f_2\Vert_{L^{p_2}(\bbr)}.
\end{equation*}
Choosing $0<\epsilon<\wt{\delta} \frac{\theta_1+\theta_2}{\theta_0}$ so that $ \wt{\delta} \theta_1+\wt{\delta} \theta_2  -\epsilon \theta_0>0$ and setting $\delta:= \wt{\delta} \theta_1+\wt{\delta} \theta_2  -\epsilon \theta_0$, we derive
\begin{equation*}
  \Big\|\sup_{\gamma\in\bbz}\big| \LL_{\Omega,\mu}^{\gamma}\big(f_1,f_2\big)\big| \Big\|_{L^{p}(\bbr) }\lesssim 2^{-\delta \mu}\Vert\Omega\Vert_{L^q(\S^{1},u_A^q)}\Vert f_1\Vert_{L^{p_1}(\R)}\Vert f_2\Vert_{L^{p_2}(\R)}. 
\end{equation*}
    
This completes the proof.

\hfill

\section{Proof of Proposition \ref{lommushweakre}}\label{Sec: Section 10}

The proof is almost same as that of Proposition \ref{PropoHLS}.
Assuming 
$$\Vert f_1\Vert_{L^1(\bbr)}=\Vert f_2\Vert_{L^1(\bbr)}=\Vert \Omega\Vert_{L^1(\mathbb{S}^1,u_1)}=1,$$
we need to prove that for any $\epsilon>0$ and for any $\lambda>0$
\begin{equation}\label{mainproponeedest}
  \bigg|\Big\{ x\in \bbr: \LL_{\Omega,\mu}^{\sharp}\big(f_1,f_2\big)(x)>\lambda\Big\} \bigg|\lesssim_{\epsilon} 2^{\epsilon \mu}\lambda^{-\frac{1}{2}},
\end{equation}
but it can be done by mimicking the proof of Proposition \ref{PropoHLS} in Section \ref{Sec: Section 7}.
To be specific, applying \eqref{CZdecompff12}, \eqref{CZdecompff12:2}, and \eqref{CZdecompff12:3}, we bound
 the left-hand side of \eqref{mainproponeedest}  by the sum of
\begin{align*}
  \mathcal{I}_1&:=\bigg|\Big\{ x\in \bbr: \LL_{\Omega,\mu}^{\sharp}\big(g_1,g_2\big)(x)>\frac{\lambda}{3}\Big\} \bigg|\\
  \mathcal{I}_2&:=\bigg|\Big\{ x\in \bbr: \sum_{Q\in\mathcal{A}}\LL_{\Omega,\mu}^{\sharp}\big(b_{1,Q},f_{2,Q}\big)(x)>\frac{\lambda}{3}\Big\} \bigg|\\
  \mathcal{I}_3&:=\bigg|\Big\{ x\in \bbr: \sum_{R\in\mathcal{B}}\LL_{\Omega,\mu}^{\sharp}\big(f_{1,R},b_{2,R}\big)(x)>\frac{\lambda}{3}\Big\} \bigg|
\end{align*}
It follows from Chebyshev's inequality and \cite[Proposition 3]{Park2025}  with $L^1(\S^1,u_1)\hookrightarrow L^1(\S^1)$ that for some $M>0$
\begin{align*}
  \mathcal{I}_1\lesssim \lambda^{-2}\big\Vert \LL_{\Omega,\mu}^{\sharp}(g_1,g_2) \big\Vert_{L^2(\bbr)}^2\lesssim_M {\mu^M} \lambda^{-2}\Vert \Omega\Vert_{L^1(\mathbb{S}^1)}^2\Vert g_1\Vert_{L^4(\bbr)}^2\Vert g_2\Vert_{L^4(\bbr)}^2\lesssim_{\epsilon,M}2^{\epsilon \mu} \lambda^{-\frac{1}{2}}.
\end{align*}
To estimate $\II_2$, we write
$$\mathcal{I}_2\le \Big| \bigcup_{Q\in\mathcal{A}}Q^*\Big|+ \bigg|\Big\{ x\in \Big(\bigcup_{Q\in\mathcal{A}}Q^*\Big)^c: \sum_{Q\in\mathcal{A}}\LL_{\Omega,\mu}^{\sharp}\big(b_{1,Q},f_{2,Q}\big)(x)>\frac{\lambda}{3}\Big\} \bigg|$$
where $Q^*$ as the concentric dilation of $Q$ with $\ell(Q^*)=10\ell(Q)$, and observe that
$$\Big| \bigcup_{Q\in\mathcal{A}}Q^*\Big|\lesssim \lambda^{-\frac{1}{2}}.$$
Moreover, since
$$\LL_{\Omega,\mu}^{\sharp}\big(b_{1,Q},f_{2,Q}\big)(x)\le \sum_{\gamma\in\bbz}\big| T_{K_{\mu}^{\gamma}}\big(b_{1,Q},f_{2,Q}\big)(x)\big|,$$
we have
\begin{align*}
  &\bigg|\Big\{ x\in \Big(\bigcup_{Q\in\mathcal{A}}Q^*\Big)^c: \sum_{Q\in\mathcal{A}}\LL_{\Omega,\mu}^{\sharp}\big(b_{1,Q},f_{2,Q}\big)(x)>\frac{\lambda}{3}\Big\} \bigg|\\
  &\lesssim\lambda^{-1}\sum_{Q\in\mathcal{A}}\int_{(\bigcup_{Q\in\mathcal{A}}Q^*)^c} \LL_{\Omega,\mu}^{\sharp}\big(b_{1,Q},f_{2,Q}\big)(x) \; dx\\
  &\le\lambda^{-1}\sum_{Q\in \mathcal{A}}\sum_{\gamma\in\bbz}\int_{(\bigcup_{Q\in\mathcal{A}}Q^*)^c}\big| T_{K_{\mu}^{\gamma}}\big(b_{1,Q},f_{2,Q}\big)(x)\big|\; dx
\end{align*}
and the last expression has been already estimated in \cite[pages 10--12]{Ho_La_Sl_submitted} by a constant depending on $\epsilon$ times
$$2^{\epsilon \mu}\lambda^{-\frac{1}{2}}.$$
This yields
$$\mathcal{I}_2\lesssim_{\epsilon} 2^{\epsilon \mu}\lambda^{-\frac{1}{2}}.$$
Similarly, we also have
$$\mathcal{I}_3\lesssim_{\epsilon} 2^{\epsilon\mu}\lambda^{-\frac{1}{2}},$$
and this completes the proof.

\hfill

\end{document}